\documentclass{amsart}

\usepackage{amssymb,amsmath,amscd,amsthm}
\usepackage{graphicx,psfrag,epsfig,multirow}

\begin{document}
\title{The Hilbert 16-th problem and an estimate for cyclicity of
an elementary polycycle}
\author{V.Kaloshin}
\thanks{The first author is partially supported by 
the Sloan Dissertation Fellowship and the American Institute 
of Mathematics Five-year Fellowship}
\address{Fine Hall, Princeton University, Princeton, NJ, 08540}
\email{kaloshin@math.princeton.edu}
\maketitle

\markboth{Bifurcation of polycycles}{V. Kaloshin}       
 
\newtheorem{Thm}{Theorem}
\newtheorem{Def}[Thm]{Definition}
\newtheorem{Lm}{Lemma}
\newtheorem{Prop}{Proposition}
\newtheorem{Rem}{Remark}
\newtheorem{Cor}{Corollary}

\def\bdef{\begin{Def}}
\def\endef{\end{Def}}
\def\bthm{\begin{Thm}}
\def\ethm{\end{Thm}} 
\def\bprop{\begin{Prop}}
\def\enprop{\end{Prop}} 
\def\blm{\begin{Lm}}
\def\elm{\end{Lm}}
\def\bcor{\begin{Cor}}
\def\ecor{\end{Cor}}
\def\brm{\begin{Rem}}
\def\erm{\end{Rem}}
\def\bfig{\begin{picture}}
\def\efig{\end{picture}}
\def\beq{\begin{eqnarray}}
\def\eneq{\end{eqnarray}}
\def\beal{\begin{aligned}}
\def\enal{\end{aligned}}
\def\om{\omega}
\def\~{\tilde}
\def\O{\Omega}
\def\geo{geometric\ }
\def\Bbb{\mathbb}
\def\R{\Bbb R}
\def\C{\Bbb C}
\def\Z{\Bbb Z}
\def\G{\Gamma}
\def\K{\mathcal K}
\def\Cal{\mathcal}
\def\cal{\mathcal}
\def\dt{\delta}
\def\Si{\Sigma}
\def\bt{\beta}
\def\eps{\epsilon}
\def\lb{\lambda}
\def\gm{\gamma}
\def\m{\mathbf m}
\def\td{\tilde}
\def\al{\alpha}
\def\asik{Khovanski}
\def\W{\Omega}
\def\land{\wedge}
\def\ell{\mathsf p} 
\def\F{\cal F}
\def\J{\bf J}
\def\j{\bf j}
\def\o{\emptyset}

\section{Introduction}

Consider a polynomial line field on the real $(x,y)$-plane
\begin{eqnarray} \label{H}
\frac {dy}{dx}=\frac {P_n(x,y)}{Q_n(x,y)},
\ P_n,Q_n - {\text { polynomials }}, \deg P_n,Q_n \leq  n.
\\
H(n)= \boxed {\text{uniform bound for the number of limit 
cycle of}\ (\ref{H}).} \nonumber
\end{eqnarray}

One way to formulate the Hilbert 16-th Problem is the following:

{\bf{Hilbert 16-th Problem (HP).}} {\it{Find an estimate for $H(n)$
for any $n \in \Z_+$.}}

We shall discuss problems related to the following:

{\bf{Existential Hilbert 16-th Problem (EHP).}}
{\it{Prove that $H(n)< \infty$ for any $n \in \Z_+$.}}

The problem about finiteness of number of limit cycles for an 
individual polynomial line field (\ref{H}) is called
{\it{Dulac problem}} since the pioneering work of Dulac
who claimed in 1923 to solve this problem, but an error
was found by Ilyashenko.

The Dulac problem was solved by two independent and rather
different proofs given almost simultaneously by 
Ilyashenko {\cite{I}} and Ecalle {\cite{E}}.
However, both proofs do not allow any generalization to 
solve Existential Hilbert Problem.

Consider the equation (\ref{H}) for different polynomials
$(P_n(x,y),Q_n(x,y))$ as the family of line fields on $\R^2$
depending on parameters of the polynomials. Using a central 
projection $\pi:\Bbb S^2 \to \R^2$ and homogenuity with respect 
to parameters of the equation (\ref{H}) (line fields 
$\lb P_n(x,y)/ \lb Q_n(x,y)$ and $P_n(x,y)/Q_n(x,y)$ for any 
$\lb \neq 0$ are the same) one can construct a {\it{finite 
parameter family of analytic line fields on the shpere $\Bbb S^2$ 
with a compact parameter base $B$ (see e.g. {\cite{IY2}} for 
details).}} After this reduction Existential Hilbert Problem
becomes a particular case of the following 

{\bf{Global Finiteness Conjecture (GFC).}}{\it{ (see e.g. 
{\cite{R}}) For any family of line fields on $\Bbb S^2$
with a compact parameter base $B$ the number of limit cycles 
is uniformly bounded over all parameter values.}}

We refer the reader to the volumes {\cite{S}} and {\cite{IY2}}
where various development of these and related problems are 
discussed. Families of analytic fields are extremely difficult
to analyze.  In the middle of 80's Arnold {\cite{AAI}} proposed to 
consider generic families of smooth vector fields on $\Bbb S^2$.
A smooth analog of Global Finiteness Conjecture is the following

{\bf{Hilbert-Arnold Problem (HAP).}} {\it{(e.g.{\cite{IY2}}) 
Prove that in a generic finite parameter of vector fields on 
the sphere $S^2$ with compact base $B$, the number of limit cycles 
is uniformly bounded. }} 

Assume for a moment that a polynomial (or a generic smooth) 
vector field on the sphere $\S^2$ has an infinite number of limit
cycles. By the Poincare-Bendixon Theorem, any limit cycle should 
surround an equilibrium point and, since our vector field has at most
finitely many equilibria, there should be an infinite ``nested'' 
sequence around one of equilibria. Then those ``nested'' limit
cycles have to accumulate (in the sense of Hausdorff metric) to
a certain contour (polygon) consisting of equilibria (as vertices) 
and separatric curves (sides of that polygon) connecting them.
Such objects are called {\it polycycles}. 
It turns out that a possible solution to Hilbert-Arnold Problem 
reduces to investigation of bifurcation of polycycles. 
Let us give several definitions.

\bdef \label{polyc}
A polycycle $\gm$ of a vector field on the sphere $\Bbb S^2$
is a cyclically ordered collection of equilibrium points 
$p_1, \dots, p_k$ (with possible repetitions) and different arcs
$\gm_1, \dots , \gm_k$ (integral curves of the vector field) 
connecting them in the specific order:
the j-th arc $\gm_j$ connects $p_j$ with $p_{j+1}$ for $j=1, \dots, k$.  
\endef

\bdef \label{cycle}
Let $\{\dot x=v(x,\eps)\}_{\eps \in B^n}, \ x \in \Bbb S^2,$ be an
$n$-parameter family of vector fields on $\Bbb S^2$
having a polycycle $\gm$ for the critical parameter value $\eps_*$.
The polycycle $\gm$ has cyclicity $\mu$ in the family 
$\{v(x,\eps)\}_{\eps \in B^n}$ if there exist neighborhoods 
$U$ and $V$ such that $\Bbb S^2 \supseteq U \supset \gm,
\ B \supseteq V \in \eps_*$ and for any $\eps\in V$ the field 
$v(\cdot,\eps)$ has no more than 
$\mu$ limit cycles inside $U$ and $\mu$ is the minimal number
with this property.
\endef

{\it{ Examples}} \ \ 
1) In a generic $n$-parameter family, the maximal multiplicity of 
a degenerate limit cycle does not exceed $n+1$, e.g. in codimension $1$ a
semistable limit cycle has cyclicity $2$. Thus, the cyclicity
of a trivial polycycle (a polycycle without singular points)
in a generic $n$-parameter family does not exceed $n+1$. 

2) (Andronow-Leontovich, 1930s; Hopf, 1940s).
A nontrivial polycycle of codimension $1$ has cyclicity at most $1$.

3) (Takens, Bogdanov, Leontovich, Mourtada, Grozovskii,
early 1970s-1993 (see {\cite{G}}, {\cite{KS}} and references there)).
A nontrivial polycycle of codimension $2$ has cyclicity at most $2$.

\bdef The bifurcation number $B(k)$ is the maximal 
cyclicity of a nontrivial polycycle occurring in a generic $k$-parameter
family.
\endef

The definition of $B(k)$ does not depend on a choice of the 
base of the family, it depends only on the number $k$ of parameters. 

{\bf{Local Hilbert-Arnold Problem (LHAP)}} {\it{e.g.{\cite{IY1}}
Prove that for any finite $k$, the  bifurcation number $B(k)$ is 
finite and find an upper estimate for $B(k)$.}}

It turns out that a solution to Local Hilbert-Arnold Problem implies 
a solution to Hilbert-Arnold Problem.

Similarly to the generic smooth vector fields, in the case of analytic 
vector fields one can define so-called {\it a limit periodic set} 
{\cite{FP}}, {\cite{R}}, {\cite{IY1}}, which is either a polycycle or 
has an arc of equilibrium points\footnote{ generic vector fields
can not have an arc of equilibrium points}, and formulate 

{\bf{Local Finiteness Conjecture (LFC)}} {\it e.g.{\cite{R}}
Prove that any limit periodic set occuring in an analytic family
of vector fields on $\S^2$ has finite cyclicity in this family}.

Smooth vector fields are more flexible then analytic vector fields
and easier to analyze. A strategy to attack Existential Hilbert
Problem proposed by Arnold {\cite{AAI}} (see also {\cite{IK}}) 
is first understand generic smooth vector fields and then 
try to apply developed methods to analytic vector fields.
Let us summarize the discussion in the form of diagramm:

\begin{figure}[htbp]
  \begin{center}
    \begin{psfrags}
     \includegraphics[width= 4in,angle=0]{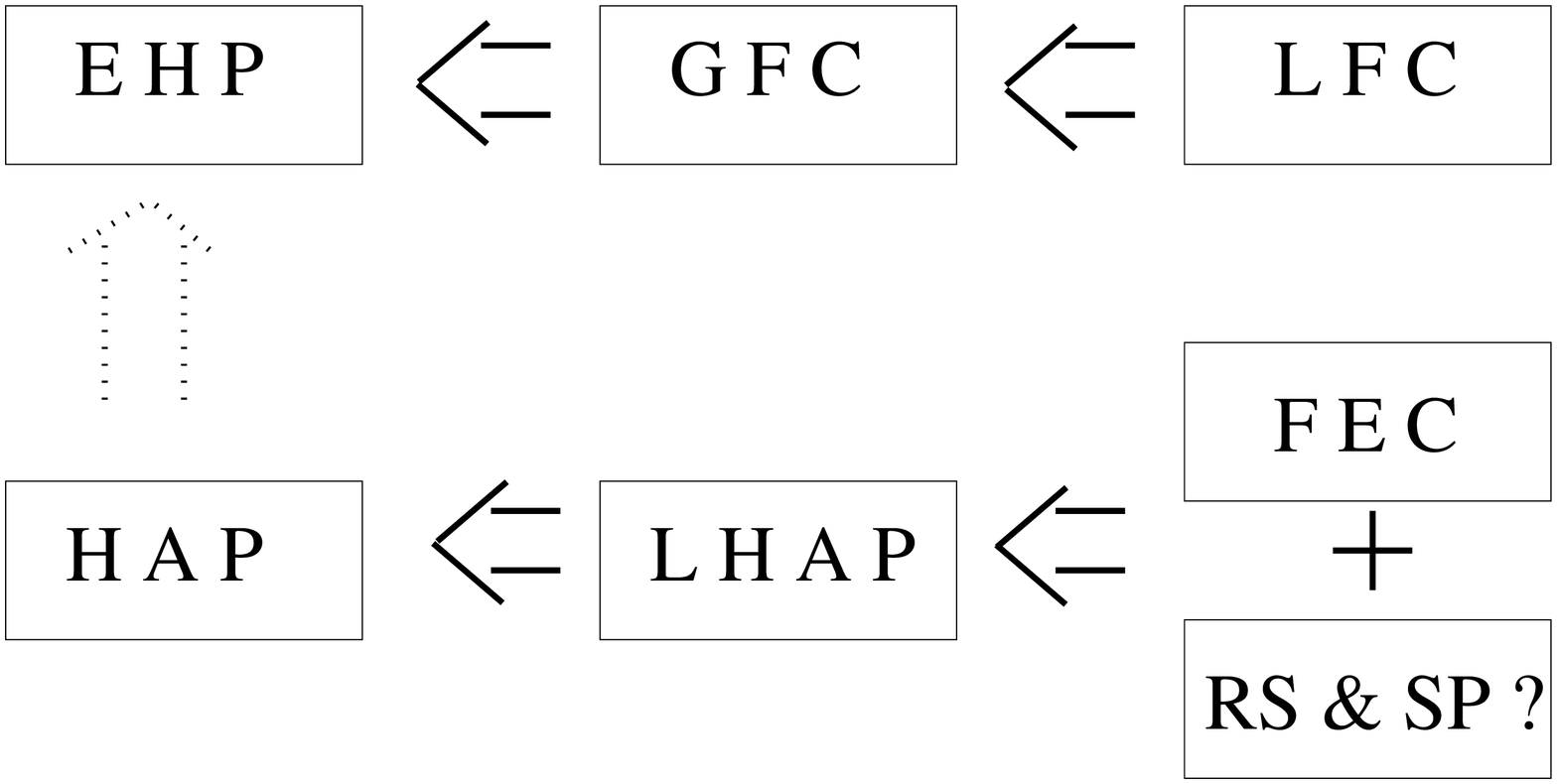}
    \end{psfrags}
    \caption{}
  \end{center}
\end{figure}

Now we shall formulate the Main Result of the paper.

\bdef
A singular (equilibrium) point of a vector field on the two-sphere
is called {\it{elementary}} if at least one eigenvalue of its
linear part is nonzero. A polycycle is called an {\it{elementary}}
polycycle if all its singularities are elementary.
\endef

The Local Hilbert-Arnold problem was solved under 
the additional assumption that a polycycle have elementary 
singularities only.

\bdef The elementary bifurcation number $E(k)$ is the maximal 
cyclicity of a nontrivial elementary polycycle occurring in a generic 
$k$-parameter family.
\endef

From examples 2) and 3) above it follows that
$$
E(1)=1, \quad E(2)=2.
$$
Information about behavior of the function $k \mapsto E(k)$ 
has been obtained recently. The First crucial step was done
by Ilyashenko and Yakovenko: 
 
{\bf{Finiteness Theorem}}\ {\it{(Ilyashenko and Yakovenko {\cite{IY3}})
For any $n$ the elementary bifurcation number $E(n)$ is finite.}}

\bcor Under the assumption that families of vector fields
have elementary singularities only the global Hilbert-Arnold conjecture 
is solved, i.e. any generic finite parameter family of
vector fields on the sphere $\Bbb S^2$ with a compact base and 
only elementary singularities 
has a uniform upper bound for the number of limit cycles.
\ecor

{\bf{Main Theorem.}}\ {\it{For any $k\in \Z_+$}} 
\beq
E(k) \leq 2^{25k^2}.
\eneq
This is the first known sufficiently general estimate for cyclicity 
of polycycle.
The case of a polycycle consisting only one singular point
with no arcs at all, is well known. An elementary equilibrium point
can generate limit cycles in its small neighborhood if it is a slow 
focus, that is the linearization matrix has a pair of 
two imaginary eigenvalues. This bifurcation was investigated by 
Takens {\cite{Ta}}. 

\bcor Under the assumption that all the polycycles are elementary
the Main Theorem gives a solution to the Local Hilbert-Arnold
problem.
\ecor

The Main Theorem is an improvement of Ilyashenko-Yakovenko Finiteness 
Theorem. It is a great pleasure for the author to say that the paper
of Ilyashenko-Yakovenko {\cite{IY3}} was a corner stone for the
present paper. In {\cite{IY3}} the authors made an extremely important 
step: they found a pass from {\it{bifurcation theory \textup{to} 
singularity theory}} using the Khovanskii reduction method
{\cite{Kh}}. We follow this pass up to some point and using some new 
ideas getting the first sufficiently general estimate for the
cyclicity of polycycles. To make this paper readable we have to reproduce
some points from {\cite{IY3}} and we are sorry for repetition,
but we think that it is necessary for a better understanding.

\subsection{Three stages of the proof}
The proof of the Main Theorem consists of three steps.
Relation to the proof of the Finiteness Theorem {\cite{IY3}}
is discussed after this short description.

{\it{Step 1. Normal forms for local families of vector fields
and their integration}}\ \ In section \ref{nform} we use
normal forms to establish an explicit form for the Poincare 
correspondence map  near equilibrium points on the polycycle
under consideration. In  {\cite{IY3}} it is shown that 
these maps satisfy Pfaffian (polynomial differential)
equations with coefficient of polynomials depending smoothly
on the parameters of the family. As the result 
a {\it{basic system}} of equations for determination of limit
cycles is obtained.

{\it{Step 2. the Khovanskii reduction method}}\ \ 
In section \ref{asik} we discuss a variation of the 
Khovanskii method  {\cite{Kh}}. This method allows us to 
investigate systems of equations that involve functions 
satisfying Pfaffian equations. In section \ref{pfreduc}
we present a formal reduction from the basic system to a {\it{mixed
functional-Pfaffian}} system which is done in {\cite{IY3}}
together with upper bounds for degrees of involved 
into the procedure polynomials. After application of the Khovanskii
method to the mixed functional-Pfaffian system we obtain several 
{\it{chain maps}}, the maps of the form
\beq \label{chain}
x \mapsto (P_1,\dots ,P_n)\circ (x,f(x),f'(x),\dots,f^{(n)}(x)),
\eneq
where $P=(P_1,\dots ,P_n)$ is a vector-polynomial given by its 
coordinate functions of known degree and $f$ is a generic function. 
The problem of estimating the number of limit cycles reduces
to estimating the number of regular preimages of some 
special points by the chain map. Special points form
an open cone-like semialgebraic set $K$ in the image.

Denote by $F$ the map
$F:x \mapsto P \circ(x,f(x),f'(x),\dots,f^{(n)}(x))$
which is called the $n$-th jet of $f$. Denote by $L_F$ 
the linearization of $F$ at point $x=0$.
 
{\it{Step 3. Bezout's theorem for the Chain maps}} 
In section \ref{apstrat}  we construct an algebraic set 
$\Sigma$ in the image of $F$ (in the space of $n$-jets).
If $F$ is transversal to $\Sigma$, then the number of 
preimages of any point $a$ from a set of special points $K$
is {\it{the same}} for $F$ and its linearization $L_F$ at
zero, namely,
\beq \label{bez}
\#\{x: P\circ F(x)=a\}=\#\{x: P\circ L_F(x)=a\}\leq \prod_{j=1}^k \deg P_j.
\eneq
But $L_F$ is a linear map and one can apply Bezout's theorem
to estimate the right-hand side of the equality.
This observation completes the proof of the Main Theorem.

Let us discuss relation of this proof to the proof of the
Finiteness Theorem by Ilyashenko \& Yakovenko {\cite{IY3}}. Step 1
of both proof is the same. We just refer to appropriate 
statements in {\cite{IY3}}. Step 2 in this proof is
slightly different for the one in {\cite{IY3}}. After application
of the Khovanskii method they obtain the same collection of chain maps
of the form (\ref{chain}). However, they investigate
the number of regular preimages of points in the image by the chain maps
{\it{without any restriction}} on those points. In the present
proof, using additional arguments in the Khovanskii method, 
we reduce consideration to
only preimages of {\it{special}} points, i.e. points from a 
tiny cone-like set in the image. At his point our proof goes 
independently, because investigation of the number of regular 
preimages of special points is more concrete problem.

Let us present a more detailed description of each step of
the proof.

\subsection{Normal forms of local families and their integration}
This step is done in {\cite{IY3}} $\S 0.3$ and $\S 1$. We just say
several words about it.

It turns out that in a small neighborhood of an elementary equilibrium 
point there exists a finitely differentiable normal coordinates
(in the Cartesian product of the phase space and the parameter space),
so-called normal forms of an equilibrium point. The list of
finitely differentiable normal forms was obtained in {\cite{IY1}}.
The main feature of the list: all normal forms are {\it{polynomial
and integrable.}} The smaller the neighborhood of a normal form,
the higher its smoothness. So smoothness can be chosen 
arbitrary large. All normal forms are summarized in Table 1
sect.\ref{nform}.

In a small neighborhood of an elementary equilibrium point
 one can choose two small segments, say $\Si^-$ and $\Si^+$,
transversal to the vector field for the critical value
of parameter and explicitly calculate the Poincare 
(correspondence) map which maps a point from one
segment say $\Si^-$ along the corresponding phase curve 
to a point from the other segment $\Si^+$ (see Fig.1).
For an appropriate choice of segments $\Si^-,\Si^+$ and 
 coordinate functions $x,$ in $\Si^-,\Si^+$  respectively, and 
a smooth function $\lb(\eps)$ in the original parameter
$\eps$ of the family the Poincare return map $\Delta_\eps:x\to y$
can be explicitly computed.
Moreover, there is a Pfaffian (with polynomial coefficients)
1-form $\om$ of the form
\beq \label{pfaff}
P(x,y,\lb(\eps))\ dx+\ Q(x,y,\lb(\eps))\ dy=0
\eneq
which vanishes on the graph $y=\Delta_\eps(x)$.
For example, in the case of a nonresonant saddle 
$\Delta_\eps(x)=x^{\lb(\eps)}$ and $\om=x\ dy+\ \lb(\eps)y\ dx$.
See Table 1 for the other cases.

\begin{figure}[htbp]
  \begin{center}
    \begin{psfrags}
      \psfrag{S3+}{$\Si_3^+$}
      \psfrag{S2+}{$\Si_2^+$}
      \psfrag{S1+}{$\Si_1^+$}
      \psfrag{S3-}{$\Si_3^-$}
      \psfrag{S1-}{$\Si_1^-$}
      \psfrag{S2-}{$\Si_2^-$}
      \psfrag{02}{$p_2$}
      \psfrag{03}{$p_3$}
      \psfrag{O1}{$p_1$}
      \psfrag{D1}{$\Delta_1$}
      \psfrag{D2}{$\Delta_2$}         
      \psfrag{D3}{$\Delta_3$}      
      \psfrag{f1}{$f_1$}
      \psfrag{f2}{$f_2$}         
      \psfrag{f3}{$f_3$}  
     \includegraphics[width= 4in,angle=0]{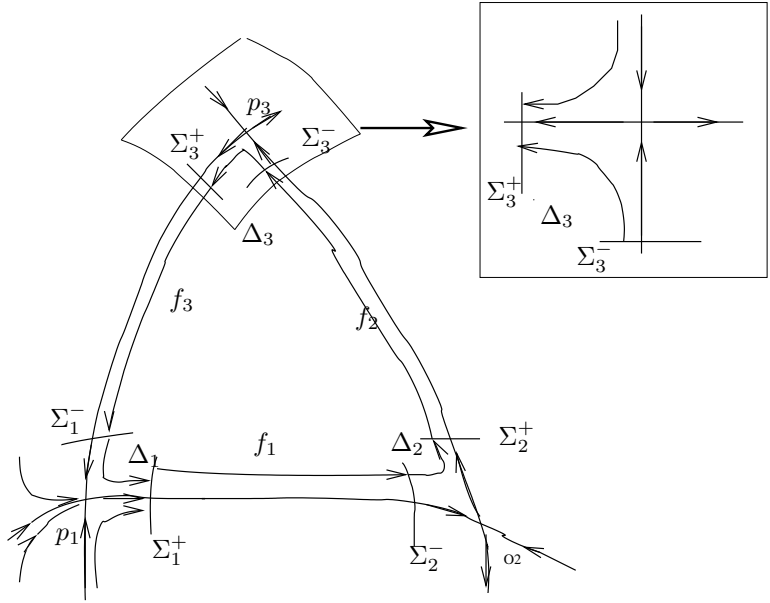}
    \end{psfrags}
    \caption{Construction of entrance and exit transversals}
  \end{center}
\end{figure}

\subsection{Singular-regular systems determining the number of
limit cycles}
We present a description of a system of equations determining the 
number of limit cycles. For a detailed description we refer
to {\cite{IY3}} $\S 0.4$ and $\S 1.4$. 

Let $\gm$ be a polycycle, occurring in a generic $k$ parameter
family, with equilibrium points $p_1,\dots ,p_n$
(possibly with repetitions) and connecting phase curves
$\gm_1,\dots ,\gm_n$ such that $\gm_j$ connects equilibria $p_j$
with $p_{j+1}$ respectively. For each 
$1\leq j \leq n$ endow the point $p_j$ with a
{\it{$C^r$-normal coordinate charts}} $U_j$. Consider 
transversal segments ``entrance'' $\Si^-_j$ and ``exit'' $\Si^+_j$
which are parallel to coordinate axis of the 
normal chart. The phase curve $\gm_{j-1}$ enters the 
neighborhood $U_j$ through $\Si^-_j$
and the phase curve $\gm_{j}$ exists $U_j$ through
$\Si^+_j$. The normal coordinates induce coordinates $x_j$ and 
$y_j$ on $\Si^-_j$ and $\Si^+_j$ respectively. For some parameter values
the corresponding vector field defines the following collection
of Poincare maps:
\beq
\begin{aligned}
\Delta_j(\cdot,\eps):x_j\to y_j=\Delta_j(x_j,\eps),\ j=1,\dots ,n\\
f_j(\cdot,\eps):y_j \to x_{j+1}=f_j(y_j,\eps), \ j=1,\dots ,n\ \ (\mod n), 
\end{aligned}
\eneq
where $\Delta_j(\cdot,\eps)$ is a local Poincare map
form the ``entrance''segment $\Si^-_j$ to the 
``exit'' segment $\Si^+_j$ and $f_j(\cdot,\eps)$ is a semilocal 
Poincare map along the phase curve $\gm_j$ form the ``exit'' segment 
$\Si^+_j$ to the ``entrance'' segment $\Si^-_{j+1}$.

Now we decompose the monodromy map (the Poincare first return map)
along the polycycle $\gm$ into the chain of the local singular maps  
$\Delta_j$ and the semilocal regular maps $f_j$ of the total length
$2n$. Limit cycles correspond to the fixed points of the monodromy.
But instead of writing one equation for the fixed points of the monodromy 
we consider a system of $2n$ equations, which will be called the 
preliminary basic system:
\beq
\begin{cases} \label{prebasic}
y_j = \Delta_j (x_j, \eps) \quad j=1, \dots, n \\
x_{j+1}=f_j(y_j,\eps), \ j=1,\dots ,n\ (\mod n)
\end{cases}
\eneq

Recall that $x_j$'s are $C^r$-normal coordinates on $\Si_j^-$ and 
$y_j$'s are  $C^r$-normal coordinates on $\Si_j^+$. Thus the system 
involves $C^r$-smooth regular functions $f_j$'s and the maps  
$\Delta_j$ from the list 
(modulo reparametrization $\eps \to \lb (\eps)$),
that are essentially singular. The problem now is to estimate the number 
of solutions uniformly over all sufficiently small parameter values.

\subsection {The Khovansky reduction method.}

The system (\ref{prebasic}) is not easy to analyze, because it has 
the singular functions $\Delta_j$. The first key idea of the second step 
is to replace these singular equations in (\ref{prebasic}) by
the Pfaffian (polynomial differential) equations in the form (\ref{pfaff}).
As a result we obtain the {\it{mixed}}  functional-Pfaffian
system of the form 
\beq
\begin{aligned} \label{funpfaf}
\begin{cases}
\om_j = 0  \\
F_j(x,y,\eps) = 0 \quad \quad \quad j=1, \dots, n \quad
\end{cases}
\\
\om_j = P_j\ dx_j + Q_j\ dy_j,\quad F_j(x,y,\eps)=x_{j+1}-f_J(y_j, \eps)\\
(x,y)=(x_1, y_1, \dots , x_n, y_n) \in (\R^{2n}, 0),\ \eps \in (\R^k, 0),
\end{aligned}
\eneq
where $\om_j$ are Pfaffian forms  in the form (\ref{pfaff}). This system 
can be interpreted as follows: one has to take an integral manifold 
$\Gamma$ for the Pfaffian equations of the system (\ref{funpfaf}) and
compute its intersection with the level set $\F^{-1}(0)$, where 
$F:\ (\R^{2n}, 0) \to \R^n$ is the map with the coordinate functions
$F_j$. In order to estimate number of isolated solutions to 
(\ref{prebasic}) one needs to estimate the number of isolated 
points in the intersection. It turns out that {\it{it is sufficient 
to analyze only transversal intersections}} of  $\Gamma$ with a generic
level set $F^{-1}(b)$ for $b$ sufficiently close to the origin in 
$\R^n$. Since the integral manifold and the level sets have complimentary 
dimensions, a transversal intersection always consists of isolated points,
which we call {\it{regular solutions}} to the system  (\ref{funpfaf}).
What we are interested in is the upper estimate for their number, 
uniform over all the integral manifolds  $\Gamma$ and all sufficiently 
small values of the parameters.

The method suggested by  A. Khovanski \cite{Kh} allows us to replace 
a mixed functional-Pfaffian system of the form (\ref{funpfaf}) by the
two systems of a similar form, but containing $n-1$ Pfaffian equations,
$n$ ``simple'' functional equations, and one {\it{special}} functional  
equation: the number of regular solutions to the initial
equation is bounded from above by the sum of the number of regular 
solutions to these two auxiliary systems.

\subsection{$a_P$-stratification and Bezout's theorem 
for a chain map $P \circ F$ with a generic $F$. }\label{subst}

In this section we shall discuss the formula (\ref{bez}). 
The problem of estimating the maximal number
of small isolated preimages is equally difficult for a chain map 
$P \circ F:\Bbb R^n \to \Bbb R^n$ with a generic map 
$F:\Bbb R^n \to \Bbb R^N$, $N\geq n$ and for a chain map 
$P \circ j^nF:\Bbb R^n \to \Bbb R^n$ with the $n$-jet of 
a generic map. We shall show that if the map $F$ 
(resp. $j^nF$) satisfies a transversality  condition 
in an appropriate space, then $F$ (resp. $j^nF$) 
can be replaced by its linear part and we can apply the Bezout 
theorem to estimate the maximal number of small inverse images
of the chain $P \circ F$ (resp. $P \circ j^nF$) uniformly over
all sequences of numbers $\eps_1, \dots, \eps_n$ decreasing
sufficiently fast to $0$. So, to simplify notations we shall 
consider a chain map of the form 
$P \circ F:\Bbb R^n \to \Bbb R^n$.

\subsubsection{A Heuristic description} \label{heur} 

Consider a chain map $P \circ F: \Bbb R^2 \to \Bbb R^2$,
where $F:\Bbb R^2 \to \Bbb R^N$ is a generic 
$C^k$ smooth map, $k>2$  and $P=(P_1,P_2):\Bbb R^N \to \Bbb R^2$ 
is a polynomial of degree $d$.
Fix a small positive $r$. We would like to estimate the maximal
number of small preimages  

\beq \label{numb}
\# \{x \in B_r(0):\ P_1 \circ F(x)=\eps,\ P_2 \circ F(x)=0\}
\eneq
for a small enough $\eps$.

To show the idea put $N=3$, $P_1(x,y,z)=x^2+y^2$, and $P_2(x,y,z)=xy$.
Assume also that $F(0)=0$. 
Denote the level set by $V_\eps=\{P_1=\eps,\ P_2=0\}$. The level set $V_\eps$
for $\eps>0$ consists of $4$ parallel lines (see Figure 2).

Notice that in our notation the number of intersections of
$F(B_r(0))$ with $V_\eps$ equals
the number of preimages of the point $(\eps,0)$
(\ref{numb}).

It is easy to see from Figure 2 that if $F$ is transversal to 
$V_0$ it is transversal to $V_\eps$ for any small $\eps>0$. Moreover,
the number of intersections $F(B_r(0))$ with $V_\eps$ equals 4
(see the points $P_1,\dots,P_4$ in Figure 2).

Another way to calculate the same number is
as follows. Let us replace $F$ by its linear part $L_F$ at zero. Then 
$\#\{x \in B_r(0):\ P_1 \circ F(x)=\eps,\ P_2 \circ F(x)=0\}=
\#\{x \in B_r(0):\ P_1 \circ L_F(x)=\eps,\ P_2 \circ L_F(x)=0\}$ and 
solving this polynomial system also yields $4$.

\begin{figure}[htbp]
  \begin{center}
    \begin{psfrags}
      \psfrag{Ve}{$V_\eps$}
      \psfrag{V0}{$V_0$}
      \psfrag{P1}{$P_1$}     
      \psfrag{P2}{$P_2$}
      \psfrag{P3}{$P_3$}           
     \includegraphics[width= 3in,angle=0]{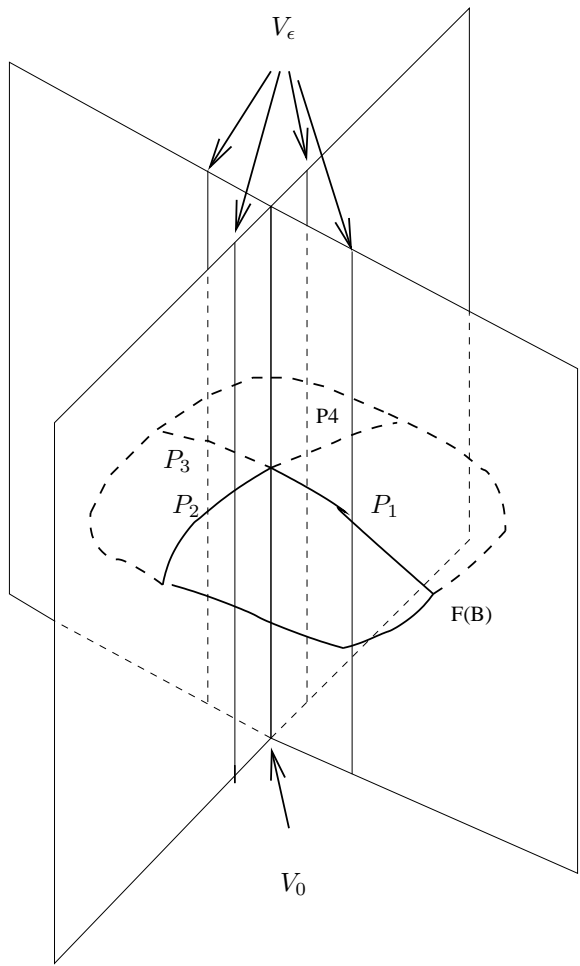}
    \end{psfrags}
    \caption{The Idealistic Example}
  \end{center}
\end{figure}

The idea behind this picture is the following:
Consider an arbitrary $N$ and a polynomial
$P=(P_1,P_2):\Bbb R^N \to \Bbb R^2$ 
of degree at most $d$, $N>2$. Define the semialgebraic variety 
$V_\eps=(P_1,P_2)^{-1}(\eps,0)$ as the level set.

Assume for simplicity that for any small $\eps \neq 0$ the level set 
$V_\eps$ is a manifold of codimension 2. 
We shall get rid of this assumption later 
(see Theorem \ref{existence} b)).
It turns out that there exists a stratification of $V_0$
by semialgebraic strata $(V_0,\Cal V_0)$ 
(a decomposition of $V_0$ into a disjoint union of semialgebraic sets 
see definition \ref{strdef}), depending on $P$ only,
such that 

\begin{multline}\label{apstra}
\boxed{F\ \text{is transversal to}\ (V_0,\Cal V_0)}
\implies \boxed{F \ \text{is transversal to}\  V_\eps}
\end{multline}
Condition (\ref{apstra}) is written  for $n=2$.
Below we shall use its analogue for an arbitrary $n$.
Let us present the key Proposition below and the simple of 
proof of it. This proof gives an insight to the main idea of the 
third step.

\bprop \label{lin} Let $B_r(a)$ be the $r$-ball centered at the point 
$a \in \Bbb R^2$ and let $L_{F,a}$ denote the linearization of $F$ at 
the point $a$. Under condition (\ref{apstra}),  
the number of intersections of the image $F(B_r(a))$ with $V_\eps$ 
coincides with the number of intersections of the 
image $L_{F,a}(B_r(a))$ with $V_\eps$, provided $r$
is small enough. That is  
\beq \label{replace}
\begin{aligned}
\#\{x \in B_r(0):\ (P_1,P_2) \circ F(x)=(\eps,0)\}=\\
\#\{x \in B_r(0):\ (P_1,P_2) \circ L_{F,a}(x)=(\eps,0)\}.
\end{aligned}
\eneq
\enprop
The argument below is independent of the codimension of $V_\eps$.
We only need condition (\ref{apstra}) and the fact that the
codimension of $V_\eps$ coincides with the
dimension of the preimage of a chain map $P \circ F$.
 
{\it{Proof}}\ \ 
Consider the $1$-parameter family of maps 
$F_t=tF+(1-t)L_F$ deforming the linear part of $F$ into $F$. 
Clearly, $F_1\equiv F$ and $F_0 \equiv L_F$. Fix a small $r>0$.
Since, $F$ is transversal to $V_0$ at $0$ all $F_t$ are 
transversal to $V_0$ at $0$. Condition (\ref{apstra})
implies that for all small $\eps$ and all $t \in [0,1]$
$F_t$ is transversal to $V_\eps$. 

Therefore, {\it{the number of intersections of $F_t(B_r(0))$ with 
$V_\eps$ is independent of $t$.}} Indeed,
assume that $\#\{F_{t_1}(B_r(0)) \cap V_\eps\} \neq  
\#\{F_{t_2}(B_r(0)) \cap V_\eps\}$ for some $t_1<t_2$.
Then as $t_1$ increases to $t_2$ there is a point $t^*$ where the
number of intersections drops or jumps.
At this point $t^*$ the condition of transversality of
$F_{t^*}$ and $V_\eps$ must fail.
This completes the proof of the proposition.

\section{Normal forms for local families and their applications.}
\label{nform}

In this section we present the functional--Pfaffian
system whose number of solutions bounds from above 
the number of limit cycles. This system was obtained in 
{\cite{IY3}}.

\subsection{ Local families and polynomial normal forms}

A local family of planar vector fields is the germ of a map,
$$
v:(\R^2,0)\times(\R^k,0)\to(\R^2,0),\qquad
(x,y,\eps)\mapsto v(x,y,\eps).
$$
A $C^r$-smooth conjugacy between two local families $v$ and $w$ of the
above form is a map
\beq \nonumber
H:(\R^2,0)\times(\R^k,0)\to(\R^2,0),\qquad
(x,y,\eps)\mapsto H(x,y,\eps),
\eneq
such that
\beq \nonumber
\ \ \ 
H_* v(x,y,\eps)=w(H(x,y,\eps),\eps),
\eneq
where $H_*$ stands for the Jacobian matrix with respect to the
variables $x,y$.  (this definition does not yet allow for
reparameterization of a local family).  Two families are finitely
differentiably equivalent, if for any $r<\infty$ there exists a
$C^r$-conjugacy between them. The two families $v,w$ are orbitally
equivalent, if there exists the  germ of a  nonvanishing function
$\phi:(\Bbb R^2,0)\times(\Bbb R^k,0)\to\Bbb R^1$ such that $v$ is
equivalent to $\phi\cdot w$.

To allow for a reparameterization of local families, we say that a
family $v(\cdot, \eps)$ is induced from another family $w(\cdot,\lb),\
\lb\in(\Bbb R^m,0)$, if $v(\cdot,\eps)=w(\cdot,\lb(\eps))$, where 
$\lb(\eps)$ is the germ of a smooth map $(\Bbb R^k,0)\to(\Bbb R^m,0)$. 
The number of new parameters $m$ may be different from $k$.

Assume that the family $w(\cdot,\lb)$ is {\it global\/} (i.e. the
expression $w(x,y,\lb)$ makes sense for all $(x,y,\lb)\in\Bbb R^{m+2}$);
this happens in particular when $w$ is {\it polynomial\/} in all its
arguments. Restricting the parameters $\lb$ onto a small neighborhood of
a certain point $(0,0,{\bf c})\in \Bbb R^2\times\Bbb R^m$, we obtain a
{\it localization\/} of the global family $w$, which 
formally becomes a local family after the parallel translation 
$\lb\mapsto\lb-\bf c$.

\bdef
1. A local family $v=v(\cdot,\lb)$ is finitely smooth orbital versal
unfolding (in short, versal unfolding) of the germ $v(\cdot,0)$, if any
other local family unfolding this germ is finitely differentiable
orbitally equivalent to a family induced from $v$.

2. A polynomial family $w(\cdot,\lb)$, $\lb\in\Bbb R^m$, is a {\it
global finitely smooth orbital versal unfolding\/} (in short, global
versal unfolding) for a certain class of local families of vector
fields, if any local family from this class is finitely differentiable
orbitally equivalent to a local family induced from some localization
of $w$. 
\endef

To investigate a versal unfolding means to investigate at the same time
all smooth local finite-parametric families which unfold the same germ
$v(\cdot,0)$.  The main result describing versal unfoldings of germs of
elementary singularities on the plane, is given by the following
theorem.

\bthm \label{nforms}{\cite{IY3}}
Suppose that a generic finite-parameter family of smooth vector fields
on the plane possesses an elementary singular point for a certain value
of the parameters. If this point has at least one hyperbolic sector,
than the family is finitely differentiable orbitally equivalent to a
family induced from some localization of one of the families given in
the second column of Table~1.
\ethm

\begin{center}
Table 1. Unfolding of elementary equilibrium points on the plane.

\end{center}
{\small
\begin{tabular}{|c |c |c |c |}
\hline
 Type &  Normal forms
 & Poincare & Pfaffian equations \\
 & & Correspondence maps & 
 \\ \hline
 & $\dot x\ =\ x,$ & & \\
 $S_0$ & $\dot y\ =\ -\lb y.$ & $y=x^\lb,$ & 
 $x\ dy\ - \lb y\ dx\ =0$ \\ &  & $x>0,\ y>0$ & \\
 & $\lb=\lb_0 \in \R^1$& & \\
 \hline 
 & $\dot x\ =\ x\ \left( \frac{n}{m}+ P_\mu(u,\lb)\right),$ & & \\
 & $\dot y\ =\ - y.$ &  & \\
 $S_\mu$& & $0\ =\ m \log y\ +$ & $y\ P_\mu(y^n,\lb)\ dx\ -$ \\
 & $u=u(x,y)=x^m\ y^n,$ & 
 $\int_{x^m}^{y^n} \frac{du}{uP_\mu(u,\lb)}.$ & 
 $\left( \frac{n}{m}+ P_\mu(y^n,\lb)\right)\times$\\
 & $P_\mu(u,\lb)=\pm u^\mu(1+\lb_\mu u^\mu)$ & $x>0,\ y>0$
  &
 $x P_\mu(x^m,\lb)\ dy=0$ \\ 
 & $+ W_{\mu-1}(u,\lb),$ &  & \\
 & $\lb=(\lb_1, \dots , \lb_\mu)$ & & 
 \\ \hline 
 & $\dot x\ =\ Q_\mu(x,\lb),$ & & \\
 & $\dot y\ =\ - y.$ & $y\ =\ C(\lb) x,$ & \\
 $D^{c}_\mu$ &  & $C={\int_{-1}^1} \frac{du}{Q_\mu(u,\lb)},$ 
 & $x\ dy\ - y\ dx=0$ \\
 & $Q_\mu(x,\lb)\ =\pm x^{\mu+1}(1+\lb_\mu x^\mu)$ & 
 $x,\ y\ \in \R^1$ & \\ \cline{1-1} \cline{3-4}
 & $+ W_{\mu-1}(x,\lb),$ & $0=\ \log y\ + $ & \\
 $D^{h}_\mu$& $\lb=(\lb_1,\dots ,\lb_\mu)$ & 
 ${\int_{x}^1}\frac{du}{Q_\mu(u,\lb)}$ &
 $Q_\mu(x,\lb)\ dy\ -$ \\
 & & $y>0,\ x\ \in \R^1$ & $y\ dx\ =0$\\ \hline
\end{tabular}}

\ 

In what follows the following notation for elementary equilibria
(the subscript indicates the degree of degeneracy):

$S_0$--- Nonresonant saddle;

$S_\mu$--- Resonant saddle whose quotient equation (the differential 
equation for $u=x^m\ y^n$ below) has the singular point of multiplicity 
$\mu+1$ at the origin, $\mu \geq 1$; if we want to specify explicitly 
the resonance between the eigenvalues, we use the extended notation 
$S_\mu^{(n:m)}$ assuming that the natural numbers $m,n$ are mutually prime;

$D_\mu$--- Degenerate saddlenode of multiplicity $\mu$;

$W_{\mu-1}(z,\lb)=\lb_0+\lb_1z+\cdots+\lb_{\mu-1}z^{\mu-1}$ is a
Weierstrass polynomial of degree $\mu-1$.

Different technical remarks concerning this table see in {\cite{IY3}}
$\S 1.1$. We just briefly describe each column.

The first two columns do not need extra words. 
In the third column of the table the Poincare correspondence maps
$y=\Delta(x,\lb)$ for the polynomial normal forms are given. They are
implicitly defined by the equations relating $x$ to $y$, these
equations depending explicitly on the parameters $\lb$ and thus
implicitly on the original parameters $\eps$. The
choice of segments transversal to the phase curves of the family
described in fig. 1.

\begin{figure}[htbp]
  \begin{center}
   \begin{psfrags}
     \psfrag{x}{\small{$x$}}
     \psfrag{y=D(x,l)}{\small{$y=\Delta(x,\lambda)$}}
     \psfrag{Dcmu}{\small{$D^c_\mu$}}
     \psfrag{Smu}{\small{$S_\mu$}}
     \psfrag{Dhmu}{\small{$D^h_\mu$}}
     \psfrag{tilq}{\small{$\tilde q$}}
     \psfrag{q}{\small{$q$}}
     \psfrag{p}{\small{$p$}}
    \includegraphics[width= 4in,angle=0]{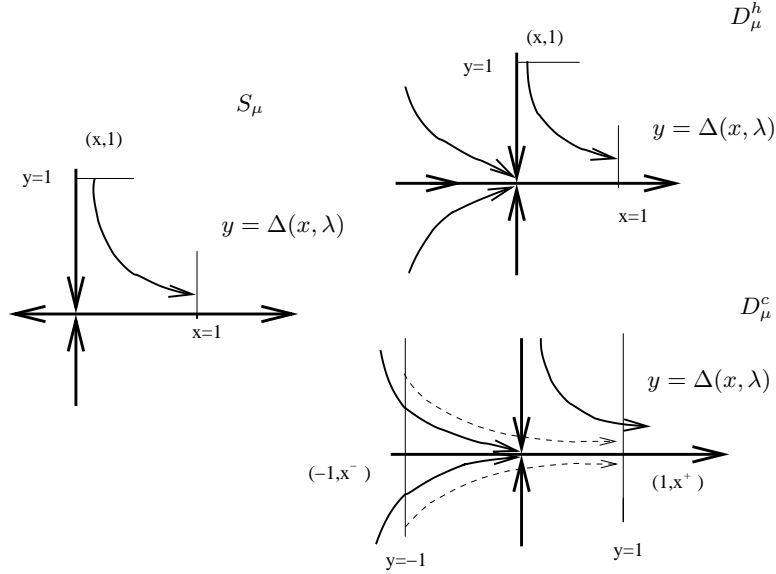}
    \end{psfrags}
   \caption{Poincare Correspondence maps}
  \end{center}
\end{figure}

\subsection{Basic system}

Here we describe the system of equations which will be analyzed from
now on. Assume that a polycycle occurs in a generic $k$-parameter
family of vector fields, and all the vertices of the polycycle are
elementary.

Then the number $n$ of vertices is $\leq k$. Moreover, one can claim that
each vertex is of one of the types $S_{\mu_j},\ \mu_j\geq 0$, or
$D_{\mu_j}^{h,c},\ \mu_j\geq 1$, and $\sum\mu_j\leq k$ 
(see {\cite{IY3}} $\S1.4$).

Next, we proceed with introducing the normalizing $C^{\mathsf p}$-smooth local
coordinates near each elementary vertex, as this is described above
(the exact order of smoothness will be specified later on).  Then a
pair of $C^{\mathsf p}$-smooth transversals may be chosen near each vertex,
and endowed with local $C^{\mathsf p}$-smooth charts $x_j,y_j$ in such a way
that the correspondence map taking a point with a coordinate $x_j$ on
the ``entrance'' transversal to a point with the coordinate $y_j$ on
the ``exit'' transversal, will be of one of the standard types listed
in Table~1.

More precisely, for each vertex $j=1,\dots,n$ Theorem \ref{nforms} 
yields the {\it localization point\/} 
${\bf c}_j=(0,\dots,0,c_j)\in\Bbb R^{\mu_j+1}$, where 
$c_j \in \Bbb R^1$ is the formal invariant of the unperturbed 
singular point, and also if $j$-th vertex is a resonant saddle, 
then the rational hyperbolicity ratio $n:m$ is explicitly specified.

Denote by $\Delta_{l,\mu}(x,\lb)$ the correspondence map for each of the
four types of singularities from Table 1, $l=S_0,S_\mu,D_\mu^c$ or
$D_\mu^h$, with the corresponding index $\mu \in \Bbb N$ (for $l=S_0$ by
definition $\mu=0$). In case $S_\mu$ with $\mu>0$ we consider the
mutually prime pair of natural numbers as an additional parameter of
the corresponding map, so in this case the rigorous notation would be
$\Delta_{S_\mu,\mu}(x,\lb;[n,m])$.

\bdef \label{f1.1}
1. The {\it unspecified basic system\/} for determination of limit
cycles occurring in $k$-parametric families of vector fields is the
system of $n$ regular and $n$ singular functional equations in $2n$
variables $x_j,y_j$, depending on parameters $\lb^j,n_j,m_j,\eps$,
\beq 
\begin{aligned}\label{1}
\begin{cases}
y_j\ \ =\Delta_{l_j,\mu_j}(x_j,\lb^j;[n_j,m_j]),\qquad
\lb^j \in\Bbb R^{\mu_j+1},
\\
x_{j+1}=f_j(y_j,\eps),\qquad \eps \in(\Bbb R^k,0).
\end{cases}
\\
j=1,\dots,n\mod(n),\quad l_j \in
\{ S_0,S_\mu,D_\mu^c,D_\mu^h \},
\\
n_j,m_j \in \Bbb N,\quad \mu_j \in \Bbb Z_+,\quad \sum \mu_j\leq k,
\quad n\leq k,
\\
\Delta\ \ {\text{depends on}}\ \ n_j,m_j \ {\text{only if}}\ \ 
l_j=S_\mu\ \ {\text{with}}\ \ \mu>0.
\end{aligned}
\eneq

2. A {\it specified\/}  basic system is one of a finite number of
unspecified basic systems together with an explicit indication of
{\it specification\/}, which by definition is the collection of:

$\bullet$ localization points ${\bf c}_j=(0,\dots,0,c_j)\in \Bbb
R^{\mu_j+1}$; in particular this means that hyperbolicity ratios of
all nonresonant saddles are explicitly given;

$\bullet$ hyperbolicity ratios $n_j:m_j$ for all resonant saddles;

$\bullet$ smooth functions $f_j(x,\eps)$ depending on the parameters $\eps$,
are defined in some open neighborhoods $(\Bbb R^{k+1},0)_j$
and $f_j(0,0)=0$;

$\bullet$ {\it characteristic size\/}, that is, the value $r>0$ which
determines the {\it domain\/} of the specified basic system as follows:
\beq
 \begin{aligned}
 (x,y)\in I_r=
 \{|x_j|<r, \ |y_j|<r,\quad j=1,\dots,n \}\subset\Bbb R^{2n};
 \\
 (\lb,\eps)\in B_r=
 \{\|\lb^j-{\bf c}_j\|<r,\ \|\eps\|<r\}\subset 
 \Bbb  R^{k+\mu_1+\cdots+\mu_n},
 \end{aligned}
\eneq
where $\lb$ is the tuple of {\it all parameters\/} of all normal forms
from Table 1, $\lb=(\lb^1,\dots,\lb^n)$; the characteristic size must be
so small that all functions $f_j$ were defined for the corresponding
values of their arguments.
\endef

{\it Notations related to definition \ref{f1.1}}\ \ 
There is only a finite number of unspecified basic systems, each one
being completely characterized by the string of discrete data
\beq \label{combtype}
\Cal T=(l_1,\mu_1,\dots,l_n,\mu_n)
\eneq
subject to the total restriction $n\le k$, $\sum\mu_j\le k$. We call
the data $\Cal T$ the {\it combinatorial type of the unspecified basic
system\/}. 

The string
\beq
\beal
\Cal S_a=(c_1,\dots,c_n,\dots,m_{j_\alpha},n_{j_\alpha},\dots)\in\Bbb
R^{n+2s}\\
 r>0,\ c_j\in\Bbb R^1,\
m_{j_\alpha},n_{j_\alpha}\in\Bbb N,\ f_j\in\Bbb C^\ell(\Bbb R^{k+1},0).
\enal
\eneq
will be referred to as the {\it algebraic part of the specification\/}
(for reasons to be clarified later), while the string of functions
$$
\Cal S_f={\bf f}=(f_1,\dots,f_n)
$$
is called the {\it functional part\/} of the specification. The functions $f_j$ are 
defined on the domain $I_r\times B_r$, where $r$ is the characteristic size
introduced earlier.

Denote by $\cal B(\Cal T,\Cal S_a, \bf f)$ the number of isolated
solutions to the specified basic system $(\Cal T,\Cal S_a, \bf f)$ in 
the domain $I_{r_0}$. One can check that $\cal B(\Cal T,\Cal S_a, \bf f)$
is defined in such a way that it bounds the cyclicity of the 
polycycle with such a specification. 

After all these notions (or rather the language) being introduced, we
may formulate the problem of estimating cyclicity of elementary
polycycles occurring in generic $k$-parametric families as follows.

\bthm \label{t1.1}
For any type $\Cal T$ of unspecified basic system and any choice of the
algebraic part $\Cal S_a$  one may
choose the order of smoothness $\ell_0$ and an open dense subset $\mathsf
F=\mathsf F_{\Cal T,\Cal S_a,r_0}$ in the space of $C^{\ell_0}$-smooth
functions $C^{\ell_0}(I_{r_0}\times B_{r_0},\Bbb R^n)$ such that for
every ${\bf f}=(f_1,\dots,f_n)\in\mathsf F$ and a sufficiently small 
characteristic 
size $r_0=r_0(\bf f)$ the number of isolated
solutions $\cal B(\Cal T,\Cal S_a,{\bf f};r_0)$ to the specified basic
system $(\Cal T,\Cal S_a, \bf f)$ in the domain $I_{r_0}$ is
uniformly bounded over all parameter values $(\lb,\eps)\in B_{r_0}$:
\beq 
\beal \label{step1}
\cal B(\Cal T,\Cal S_a,{\bf f};r_0)=\sup_{(\eps,\lb)\in B_{r_0}}
\#\{(x,y)\ {\text{satisfying (\ref{1})}}, \ (x,y)\in
I_{r_0}\}< 2^{25k^2}
\enal
\eneq
and, therefore, $E(k) \leq 2^{25k^2}$. 
\ethm

\section{The \asik\ reduction method.}\label{asik}

In this section we describe the method of reducing a
functional--Pfaffian system to a chain map
of the form (\ref{chain}).  The construction in
its full generality is described in the book \cite{Kh}.
Our exposition relies on the one in  {\cite{IY3}}, but has 
new important features so we can't just refer to neither 
\cite{Kh}, nor {\cite{IY3}}.

\subsection{Pfaffian systems and their separating solutions}

Let $M$ be a smooth orientable $n$-dimensional manifold, not
necessarily compact or connected, and $\om$ be a 
 smooth 1-form on it.

\bdef
A codimension 1 smooth submanifold $\G\subset M$ is the separating
solution for the Pfaffian equation $\om=0$, if:

a) $\G$ is the integral manifold, that is, the restriction of $\om$
on the tangent bundle of $\G$ is identically zero:
\beq \nonumber
\forall x\in\G,\ \forall v\in T_x\G\quad \om(v)=0;
\eneq

b) $\G$ does not pass through singular points of $\om$:
\beq \nonumber
\forall x\in \G,\ \exists v \in T_x M\quad \om(v)|_{T_x M}\ne0;
\eneq

c) $\G$ is the boundary of a domain $D\subseteq  M$ and the coorientation
induced on $\G$ by $\om$, coincides with its coorientation as the
boundary. In other words, on any vector pointing outward from
$D$, the form is positive.
\endef

Let now $\om_1,\dots,\om_k$ be an ordered $k$-tuple of smooth 1-forms on
$M$. Consider the system of Pfaffian equations
\beq\label{2.1}
\om_1=0,\quad \dots \quad , \om_k=0.
\eneq

\bdef
A submanifold $\G$ is the separating solution for the system of
Pfaffian equations, if there exists an increasing
chain of smooth submanifolds,

\beq\label{2.2}
\G=\G_k\subset\G_{k-1}\subset\cdots\subset\G_1\subset\G_0=M
\eneq
such that for any $j=1,\dots,k$ submanifold $\G_j$ is the separating 
solution for the Pfaffian equation on $\G_{j-1}$, determined by the 
restriction of the form $\om_j$ on the latter submanifold.
\endef

Let $\F:M \to \Bbb R^s$ be a  smooth map  $s<n-k$. 
Recall that a point $a \in {\Bbb R}^s$ is called {\it a regular value} 
for the map $\F$ if the linearization matrix, denoted by $J_\F(x)$, has 
full rank for any $x \in \F^{-1}(y)$. By the rank theorem the level 
set $V_a=\F^{-1}(a)$ of a regular value is a smooth manifold of 
dimension $n-s$.

We call $a\in {\Bbb R}^s$ {\it a regular value for $\F$ 
with respect to  Pfaffian equations (\ref{2.1})} if
$a$ is a regular value of $F$ and the $k$-form 
$\W=\om_1 \wedge \dots \wedge \om_k$, restricted to $V_a$
$\W|_{V_a}$ is nondegenerate, i.e.,  singular points of 
$\W|_{V_a}$ have measure zero.

Consider a pair of smooth maps 
$\F:M \to {\Bbb R}^s$  and $F:M\to \Bbb R^{n-k-s}$.
Now we add to a Pfaffian system (\ref{2.1}) two types of functional 
equations. The first type consists of functional 
equations $\F=a$, where $a \in \Bbb R^s$ is {\it a fixed regular value of 
$\F$ with respect to a Pfaffian system (\ref{2.1})}.
The second type consists of functional 
equations $F=b$, where $b \in \Bbb R^{n-k-s}$ is {\it a variable}.
We call equations $\F=a$, with a fixed $a \in {\Bbb R}^s$, by {\it rigid} 
equations and $F=b$, with a varying $b \in {\Bbb R}^{n-s-1}$, 
by {\it loose} equations.

\bdef
Let $\W=(\om_1,\dots,\om_k)\in(\Lambda^1(M))^k$ be a $k$-tuple of
smooth 1-forms, $\F:M \to {\Bbb R}^s$ and $F:M\to {\Bbb R}^{n-k-s}$ be
smooth maps, and $a \in {\Bbb R}^s$ be a regular value for $\F$
with respect to the $k$-tuple of smooth 1-forms. 
A solution to the {\it mixed functional--Pfaffian system\/}
\beq\label{2.3}
\W=0,\quad F=b,\quad \F=a,\qquad b\in \Bbb R^{n-k-s}
\eneq
is a pair $(\G^a, L_{b})$, where $L_{b}\subseteq  M$ is the preimage
$F^{-1}(b)$ and $\G^a$ is a separating solution for the Pfaffian 
system $\W=0$, restricted to $V_a$, and the intersection 
$\G^a \cap L_b$ is nonempty. 
\endef

The solution is {\it regular\/}, if $\G^a$ is the separating 
solution for the restriction of Pfaffian equations to $V_a$
and $b$ is the regular value for the restriction of the
map $G$ on $\G^a$.
If $(\G^a,L_b))$ is a regular solution, then the intersection 
$\G^a \cap L_{b}$ is transversal and consists of isolated points.

\bdef
The {\it \asik\ number\/} $\K\{\W,F;\F=a\}$ for the mixed system (\ref{2.3}) 
is the upper bound for the cardinalities $\#\{\G^a\cap L_{b}\}$ over all 
regular solutions of the system.
\endef

{\bf Remarks}
{\it
1. The \asik\ number is also defined if $k=0$ (resp. $s=0$), i.e., 
there are no Pfaffian (resp. rigid) equations at all. In this case 
one may put formally $\G=M$ (resp. $V_a=M$), and $\K\{\om,F;\F=a\}$ 
(resp. $\K\{\W,F;\o \}$) is equal to the upper bound
of the cardinality of preimages $\#\{L_b \cap V_a\}$
(resp. $\#\{L_b \cap \G\}$) of {\it regular values\/} for the map 
$G|_{V_a}:V_a\to {\Bbb R}^{n-k-s}$.

2. If we want to stress in the notation the phase space $M$ of the
functional--Pfaffian system, we use the notation $\Cal K_M\{\W,F,\F=a\}$.
Usually this is necessary when $F$,\ $\F$, and $\W$ are defined on the
Euclidean space ${\Bbb R}^n$, while we are interested only in solutions
belonging to some (open) ball. 

3. If we fix a coordinate system in $\Bbb R^{n-k-s}$, denote 
by $F_1, \dots ,F_s$ coordinate functions of the map 
$F:M \to \Bbb R^{n-k-s}$, and introduce the $(n-k-s)$-tuple of $1$-forms
$\W_F=(dF_1, \dots ,dF_s)$, then we can consider the following mixed system 
\beq \label{redpff}
\W=0,\quad \W_F=0, \quad \F=a.
\eneq
Regularity in the definition of the \asik\ number $\K\{\W,F;\F=a\}$ implies that
$\K\{\W,F;\F=a\}=\K\{(\W,\W_F),\o;\F=a\}$.}
 
The goal is using the \asik\  reduction principle 
estimate the \asik\ number for the mixed functional-Pfaffian
system by a linear combination of the \asik\ number for some number of 
entirely rigid functional systems.

The first step of the reduction principle is to estimate the \asik\
number for a given mixed system by a linear combination of the \asik\ numbers of two auxiliary
systems {\it containing a reduced by one number of Pfaffian equations and
an increased by one number of rigid equations\/}. 

The second step is using remark 3 replace all loose functional equations 
for pfaffian equations and apply the reduction principle to the mixed system 
consisting of $(n-s)$ Pfaffian equations $(\W,\W_F)$ and $s$ rigid equations.
Thus, after $(n-s)$ steps of the reduction principle we obtain 
a finite collection of entirely rigid functional systems.

\subsection{The Reduction principle for one Pfaffian equation}

We show how to eliminate the Pfaffian equation from the mixed system
with $(n-s-1)$ loose equations and $s$ rigid functional equations. 
\beq\label{2.4}
\om=0,\quad F=b,\quad \F=a,\qquad F:M\to \Bbb R^{n-s-1},
\qquad \F:M\to \Bbb R^{s},
\eneq
We shall outline only the key ideas.

\bdef
A smooth positive function $\rho:M\to \Bbb R_+$ is called {\it covering\/},
if it tends to zero along any nonaccumulating sequence of points in
$M$. In other terms, $\rho$ vanishes ``at infinity" on $M$, so that all
level hypersurfaces of the covering function are compact subsets of
$M$.
\endef

\brm
This definition applies both to compact and noncompact manifolds, but
in the compact case  a smooth function is covering if and only if it is
everywhere positive, thus automatically bounded away from zero.
\erm

Suppose that the manifold $M$ is endowed with the Riemann volume. Since
it is orientable, one may use the duality between functions and
$n$-forms on $M$.  Denote by the asterisk the operator taking an
$n$-form into the function (dividing by the volume form).

Fix Euclidean structures in $\Bbb R^{n-s-1}$ and $\Bbb R^s$.
Let $F_1,\dots,F_{n-s-1}$ and $\F_1,\dots,\F_{s}$ be the coordinate 
functions of the maps $F$ and $\F$ in (\ref{2.4}) respectively. 

\bdef
The contact function for the mixed
system (\ref{2.4}) is
\beq\label{2.5}
\F_{s+1}=*(\om\land dF_1\land\cdots\land dF_{n-s-1}
\land d\F_1\land\cdots\land d\F_{s}).
\eneq
The operator taking the mixed system $(\om,F;\F)$ into the 
corresponding contact function, will be denoted by
$\sigma:(\om,F;\F)\mapsto \sigma(\om,F;\F)=\F_{s+1}$.
\endef

Define the two maps by their coordinate functions,
\beq
\F^c=(\F_1,\dots,\F_{s},\F_{s+1}),\qquad
\F^\infty=(\F_1,\dots,\F_{s},\rho),
\eneq
both taking $M$ to $\Bbb R^{s+1}$, where $\F_{s+1}$ is the contact 
function (\ref{2.5}), and $\rho$ is the covering function.

\bthm \label{Kh1}
Suppose that the system (\ref{2.4}) admits regular solutions in the
sense of Definition \ref{2.3}. Then for any 
sufficiently small regular $\eps$
\beq\label{2.6}
\K\{\om,F;\F=a\}\leq \frac12 \K\{\om,F;\F^\infty=(a,\eps)\}+
\K\{\om,F;\F^c=(a,\eps)\},
\eneq
where regularity of $\eps$ means that $(a,\eps)$ is 
a regular value for both $\F^\infty$ and $\F^c$ and is 
necessary to the right-hand side systems being well defined.
\ethm

Before proving this theorem recall the Rolle lemma from an
elementary calculus.

\blm \label{Rolle} 
Consider $C^2$ Morse functions $f:S^1 \to  \Bbb R^1$ on the circle 
and $g:[0,1] \to  \Bbb R^1$ on the segment, i.e., functions $f$ and $g$
have only nondegenerate critical points. Then  for all
$a \in \Bbb R$
\beq\label{Rol}
\beal
 \#\{ x:\ f(x)=a\} \leq \#\{ x:\ f'(x)=\eps \} \\
\#\{ x:\ g(x)=a\} \leq \#\{ x:\ g'(x)=\eps \}+1
\enal
\eneq
for any sufficiently small $\eps$.
\elm
{\it{Proof}}\ \  Prove the formula for $f:S^1\to\R^1$. 
For a sufficiently small $\eps$ the number of local
maxima and minima equals $\#\{ x:\ f'(x)=\eps \}$.
Between any two consecutive preimages $x_1$ and $x_2$ 
of a point $a$, i.e., $f(x_1)=f(x_2)=a$
there exists a local minimum or maximum. Q.E.D.

Formula (\ref{Rol}) in the case of one equation 
transfers a loose equation into a rigid one.

{\it{Proof of theorem \ref{Kh1}}}\ \ 
Take a regular solution $(\G^a,L_b)$ for (\ref{2.4}), where
$L_b=F^{-1} (b)$, and suppose that the intersection $\G^a \cap L_b$ 
consists of isolated, say $d$, points. Since, $b$ is regular value 
of the restriction $F|_{\G^a}$, any small variation of $b$ may only 
increase the number of intersections. Take $b$ to be a
regular value of the restriction $F|_{V_a}$ or equivalently
$(b,a)$ to be a regular value of the map $(F,\F)$ (rather than
of the restriction of $\F$ to $\G^a$).

Then any level set $L_b$ is a one dimensional smooth manifold,
intersecting $\G$ transversally. By the classification theorem 
for one-dimensional manifolds, $L_b$ is the union of
compact (circles) and noncompact (lines) components.
Fix some orientation on each circle and each curve in $L_b$.
Consider the function $f':L_b \to \Bbb R$ which maps a point 
$x \in L_a$ to the value of the $1$-form $\om$ on the unit 
positively oriented vector tangent to $L_b$ at point $x$.

Fix a connected component, denoted by $\gm \subset L_b$.
Between any two consecutive intersection $x$ and $y$ of
$L_a$ with $\G$ values $f'(x)$ and $f'(y)$ must have different 
signs. Now we can apply the Rolle lemma with
$f'=f'$, when $\gm$ is a circle, and $f'=g'$, when $\gm$ is
a line.

Each point $x$ where $f'(x)=0$ (resp. $f'(x)$ is small) is the
point where the linear functionals 
$dF_1(x),\dots, dF_{n-s-k}(x),d\F_1(x), \dots, d\F_s(x)$. and
$\om(x)$ are linear dependent (resp. almost dependent), i.e. 
$\F_{s+1}(x)=0$ (resp. $\F_{s+1}(x)=\eps$).
This completes the proof of the theorem. Q.E.D.

\bcor
If the manifold $M$ is compact, then 
for any sufficiently small regular $\eps$
\beq\label{2.7}
\K\{\om,F;\F=a\}\leq  \K\{\om,F;\F^c=(a,\eps)\},
\eneq
where regularity of $\eps$ means that $(a,\eps)$ is 
a regular value for $\F^c$.
\ecor

{\it{Proof}}\ Indeed, in this case the first term in (\ref{2.4}) disappears.

\brm
The choice of the Riemann volume form is not essential for the above
construction. Indeed, if the volume form ${\text vol}^n$ is
replaced by a new one $b\cdot {\text vol}^n$, where $b$ is a
positive function, then the function $\F_{s+1}$ will be replaced by
$\tilde \F_{s+1}=b^{-1} \F_{s+1}$, and the map 
$\tilde \F^c=(\F_1,\dots,\F_{s},\tilde \F_{s+1})$ will have the same zero set.
\erm

\subsection{The \asik\ reduction in the general case}

Consider now the general case of the mixed system (\ref{2.3}) with $k>1$.
Suppose that $\G^a$ is a separating solution for the Pfaffian system
$\W=0$ restricted to $V_a$. By definition, this means that there exists a separating
solution $\G_k^a \subset V_a$ to the Pfaffian equation $\om_k=0$ on a separating
solution $\G_{k-1}^a \subset V_a$ to the Pfaffian system $\W'=0$ 
restricted to $V_a$, where
$\W'=(\om_1,\dots,\om_{k-1})$. Note that if $\rho$ is a covering function
on the manifold $M$, then its restriction on $\G_{k-1}^a$ is the covering
function for the latter submanifold. Next, one can endow $V_a$ 
(resp. $\G^a_{k-1}$) by the Riemann $(n-s)$-volume 
(resp. $(n-k-s+1)$-volume) form ${\text vol}^{n-s}_{V_a}$
(resp. ${\text vol}^{n-k-s+1}_{\G^a_{k-1}}$) in such a way that

\beq
\begin{aligned} \, \, \, \, \,
d\F_1\land\cdots\land d\F_{s}\land\operatorname{vol}^{n-s}_{V_a}=
\operatorname{vol}^{n}_M \ \ \ 
\om_1\land\cdots\land\om_{k-1}\land\operatorname{vol}^{n-k+1}_{\G^a_{k-1}}=
\operatorname{vol}^{n-s}_{V_a}
\end{aligned}
\eneq

Since the forms $\om_j,\ j=1,\dots,k-1$ are linear independent in a
neighborhood of $\G^a_{k-1}$, these formulas define volume forms near
$V_a$ and $\G^a_{k-1}$ respectively.  As this was mentioned before, 
the choice of the Riemann volume form does not affect the assertion 
of Theorem \ref{2.1}.

Thus one can apply Theorem \ref{2.1} to the mixed system

\beq
\om_k=0,\quad F=b, \quad \F=a 
\eneq
on the manifold $\G^a_{n-k}\subset V_a$. To describe the result, we introduce the
following two maps from $M$ to $\Bbb R^{n-k+1}$,

\beq\label{2.8}
\F^c=(\F_1,\dots,\F_{s},\rho),\qquad
\F^\infty=(\F_1,\dots,\F_{s},\F_*),
\eneq
where $\rho$ is the covering function on the manifold $M$, and
$F_*:M\to \Bbb R$ is the smooth function obtained as

\beq
\begin{aligned}\label{2.9} 
\quad \F_*=\sigma(\W,F;\F)=*(dF_1\land \dots \land
dF_{n-k-s}\land \\
\land d\F_1\land \dots \land
d\F_{s}\land \om_1\land \dots \land\om_k).
\end{aligned}
\eneq

The above choice of the Riemann volume on $\G^a_{k-1}$ implies that the
asterisk operator in the ambient manifold $M$ agrees with the asterisk
operator relevant to $\G^a_{k-1}$, therefore the formula (\ref{2.9}) defines
the same function as the formula (\ref{2.5}):  $\om_k\land
dF_1\land\cdots\land dF_{n-k-s}\land  d\F_1\land\cdots\land
d\F_{s}=\F_* \cdot\operatorname{vol}^{n-k+1}_{\G^a_{k-1}}$.

\bthm \label{Khmulti}
Let $\W,F,\F^c$, and $\F^\infty$ be as above. Then for 
any sufficiently small  regular $\eps$

\beq\label{2.10}
\K\{\W,F;\F=a\}\leq \tfrac12\K\{\W',F;\F^\infty=(a,\eps)\}+
\K\{\W',F;\F^c=(a,\eps)\},
\eneq
where regularity of $\eps$ means that $(a,\eps)$ is a regular 
value for both $\F^\infty$ and $\F^c$.
\ethm

\bcor
If either $V_a$ is compact or the restriction  
$F|_{V_a}:V_a \to \Bbb R^{n-k-s}$ is a proper map, 
i.e. preimage of any point is compact, then for any sufficiently small 
regular $\eps$
\beq\label{2.10c}
\K\{\W,F;\F=a\}\leq \K\{\W',F;\F^c=(a,\eps)\},
\eneq
where regularity of $\eps$ means that $(a,\eps)$ is a regular 
value for $\F^c$.
\ecor

{\it{Proof}}\ \ 
Straightforward application of Theorem \ref{Kh1}.

Iterating the above two statements, one can replace one by one the
Pfaffian equations by the rigid functional ones, obtaining new systems whose
\asik\ numbers estimate from above that of the initial one, by virtue
of the inequalities (\ref{2.4}) and its compact counterpart (\ref{2.6}). 
On each step one has two possibilities, either to replace a Pfaffian 
equation by the contact function, or by the covering function.  But once 
the covering function appears among the rigid functional equations, 
the level sets $F^{-1}(\cdot) \cap V_a$ becomes compact as a submanifold of
a compact $V_a$, hence on the next steps the Corollary to Theorem \ref{Khmulti} 
applies rather than the Theorem itself.

Denote by \ $T^c_\eps$ and\  $T^\infty_\eps$ the two operators, transforming the
mixed system $\{\W,F;\F=a\}$ into the mixed systems $\{\W',F;\F^c=(a,\eps)\}$ and
$\{\W',F,\F^\infty=(a,\eps)\}$ respectively, where the maps $\F^c$ and 
$\F^\infty$ are given by (\ref{2.8}) and (\ref{2.9}):
\beq
\begin{aligned}\label{2.11}
T^c_\eps\{\W,F;\F=a\}=\{\W',F;(\F,\sigma\{\W,F,\F\})=(a,\eps))\}, \\
T^\infty_\eps\{\W,F;\F=a)\}=\{\W',F;(\F,\rho)=(a,\eps))\}.
\end{aligned}
\eneq

If we start with the mixed functional--Pfaffian system $\{(\W,\W_F),\om;\F=a\}$, with
$(\W,\W_F)$ being an $(n-s)$-tuple $(\om_1,\dots,\om_k,dF_1,\dots ,dF_{n-k-s})$, 
and eliminate subsequently the forms 
$\om_k$, $\om_{k-1}, \dots, \om_1,$ $dF_1,\dots ,dF_{n-k-s}$, then the 
following maps from $M$ to $\Bbb R^n$ arise:

a) the map $\F_{[0]}$, if on each step the contact function was used,
\beq\label{2.12}
\{\om,\F_{[0]}=(a,\eps^{n-s}))\}=
(T^c_{\eps_{n-s}}\circ \*\cdots\* \circ T^c_{\eps_1})\{\W,F;\F=a\},
\eneq
where $\eps^{n-s}=(\eps_1,\dots,\eps_{n-s})$;

b) the maps $\F_{[j]}$, if on the $j$th step the covering function
was used, while on all other steps the contact ones were,
$j=1,\dots,n-s$,
\beq
\begin{aligned}\label{2.13} \qquad 
\{\om; \F_{[j]}=(a,\eps^{n-s}))\}=
(T^c_{\eps_{n-s}} \circ \cdots \circ T^c_{\eps_{j+1}}&
\circ \*T^\infty_{\eps_j} \circ \\
\circ \* T^c_{\eps_{j-1}}\circ &\cdots \circ T^c_{\eps_1})\{\W,F;\F=a\}.
\end{aligned}
\eneq

Then inductive application of Theorem \ref{Kh1} 
immediately yields the following fundamental result.

\bthm \label{mkh}
The \asik\ number for the mixed system (\ref{2.3}) on a manifold $M$
with the covering function $\rho$ and any sufficiently fast 
decaying to zero sequence $(\eps_1,\dots,\eps_{n-s})$ admits the upper estimate 
by a linear combination of \asik\ numbers of some $(n-s+1)$ auxiliary systems, 
each of them containing only rigid equations and no Pfaffian equations at all:

\beq\label{2.14} \qquad 
\K\{\W,F;\F=a\}\leq \K\{\om;\F_{[0]}=(a,\eps^{n-s})\}+
\frac12\sum_{j=1}^{k}\K\{\om;\F_{[j]}=(a,\eps^{n-s})\},
\eneq
where the maps $\F_{[j]}$ are defined by the formula
(\ref{2.11})-(\ref{2.13}).
\ethm

\subsection{Applications}\label{applications}

The  \asik \ reduction process is constructive. This leads to the 
result, which will be now formulated.

Assume that the manifold $M$ is an open domain  in $\Bbb R^n$ and admits a
polynomial covering function $\rho$. The main example is the unit
ball $B=\big\{x\in\Bbb R^n : \sum_j x_j^2<1\big\}$, for which one may take
$\rho(x)=1-\sum_j x_j^2$. Then the Riemann volume form can be chosen
algebraic, $dx_1\land\cdots\land dx_n$.

Assume also that all the forms $\om_i,\ i=1,\dots,k$ are polynomial
(i.e.~with polynomial coefficients), and the maps $\F$ and $G$ are at least
$C^{n-s}$-smooth. Then, since the operators $T^c_\eps$ and $T^\infty_\eps$ 
introduced above, involve only algebraic operations and differentiation of
functions, the following holds.

\bthm \label{appl}
If the system (\ref{2.3}) has no rigid functional equations at all ($s=0$) and 
is defined on a semialgebraic subset $M\subseteq \Bbb R^n$, 
all Pfaffian forms and the covering function
$\rho$ are polynomial of degrees $\leq d$, then all the maps
$\F_{[\alpha]}:M\to \Bbb R^n$, $\alpha=0,1,\dots,k$ constructed in Theorem
\ref{2.3} are of the form

\beq\label{2.15}
\F_{[\alpha]}=P_\alpha \circ j^{n-s} F,
\eneq
where $j^{n-s} F$ is the $(n-s)$-jet extension of $F$, 
and $P_\alpha$ are certain polynomials defined on the jet space
$J^{n-s}(\Bbb R^n,\Bbb R^{n-k})$.
For all $\alpha=0,1,\dots,k$ the degrees of the polynomial $P_\alpha$ 
admits the upper estimate by $2^\alpha (dk+n)$ and each map 
$\F_{[\alpha]}$ has a regular point for a generic map $F$.
\ethm

{\it{Proof}}\ 
The reduction procedure of elimination of a Pfaffian equation 
boils down to consecutive  application $(n-s)$ times of one of the operators
$T^c_\eps$ or $T^\infty_\eps$ (\ref{2.11}).
Consider the first step.
\beq 
\begin{aligned} \nonumber
\F^{1,\alpha}_*=\sigma\{(\W,\W_F), \o;\F)\}=*(\om_1\land\cdots\land\om_k\land 
dF_1\land\cdots\land
dF_{n-k}&)=P^{1,\alpha} \circ j^1F,
\end{aligned}
\eneq
where $P^{1,\alpha}:J^1(\Bbb R^n, \Bbb R^{n-k}) \to \Bbb R$ is a polynomial 
of degree at most $dk+n$ and is defined on the space of $1$-jets
$J^1(\Bbb R^n, \Bbb R^{n-k})$.

Denote by $\W^*=(\W,\W_F)$ the $(n-s)$-tuple of the 1-form,
$dF_s$ by $\om_{k+s}$ for $s=1,\dots,n-s$, and
the $(n-s-r)$-tuple of the 1-form, which consists of
all of 1-forms of $\W^*$ except of the first $r$, by $\W^*_r$.
Consider $\F^{r,\alpha}_*=\sigma\{\W^*_{r-1},\om_r;\F^{r-1,\alpha}\}$
and $\F^{r,\alpha}=(\F^{r-1,\alpha},\F^{r,\alpha}_*)$ for 
$r=1,\dots,n-s$. It is easy to see that $\F^{r,\alpha}_*$ has the form
$\F^{r,\alpha}_*=P^{r,\alpha} \circ j^rF$.

Using induction in $r$ it is easy to see that for $r \neq \alpha$ 
the degrees of corresponding polynomials  $P^{r-1,\alpha}$ and $P^{r,\alpha}$,
defined above, satisfy the following inequality
$\deg P^{r,\alpha} \leq 2 \deg P^{r-1,\alpha}$.
For $r=\alpha$ the operator $T^\infty_\eps$ will not exceed  
the degree $dk+n$ and $\F^{r,\alpha}=\rho$.
This implies that $\deg P^{r,\alpha} \leq 2^\alpha (dk+n)$ and complete
the proof.

\section{Functional-Pfaffian system for limit cycles}
\label{pfreduc}

In this section we consider a specified basic system $(\Cal T,\Cal
S_a,\bf f;r)$ obtained from the unspecified basic system (\ref{f1.1}), that
is we consider a system (\ref{f1.1}) together with a collection of formal
invariants $(c_1,\dots,c_n)$ of all singularities (which determines a
point in the $\lb$-space), a collection of hyperbolicity ratios
$n_{j_\alpha}:m_{j_\alpha}$ of all resonant saddles and a tuple of
sufficiently smooth functions $f_j$, on a sufficiently small open cube
$I_r\times B_r$ in the $(\eps,\lb)$-space.

Our local goal is to reduce this system to a
functional--Pfaffian system having the form described in section
\ref{asik}, with the
following properties:

$\bullet$ the new system has the form allowing for application of
Theorem \ref{Khmulti};

$\bullet$ the number of {\it regular solutions\/} to the
functional--Pfaffian system is greater or equal
to the number of {\it isolated\/} solutions to ({\ref{f1.1}), up to 
$k$, where $k$ is the number of parameters of the original family.

After application of Theorem \ref{Khmulti} we will obtain a number of 
{\it chain maps\/} with controlled degrees of the exterior polynomial parts.

\subsection{Upper estimate of the number of solutions for the
basic system: statement of results}
First of all we make the following remark. The algebraic part of the
specialization can be identified with a point
\beq \label{3.0}
\Cal S_a=(c_1,\dots,c_n,n_{j_1},m_{j_1},\dots,n_{j_s},m_{j_s})\in\Bbb
R^{n+2s},
\eneq
where $s\leq n$ is the number of resonant saddles on the polycycle: the
fact that the numbers $n_{j_\alpha},m_{j_\alpha}$ are in fact natural
will become inessential for our constructions.

\bthm \label{b->f-pf}
{\text (reduction from basic to functional--Pfaffian system)}
Consider an unspecified basic system \ref{f1.1} of a certain type $\Cal
T$ in codimension $k$, together with an arbitrary specification
$$
\Cal S=(S_a, {\bf f}; r),\quad \Cal S_a \in \Bbb R^{n+2s},\quad
{\bf f} \in C^\ell(I_r\times B_r,\Bbb R^n),\quad \ r>0.
$$

Then one can explicitly construct a functional--Pfaffian system of the
form $\{\W,F\}$, $\W=(\W_1,\dots,\W_{n+2s})$, $F=(F_1,\dots,F_{n+k+m})$,
$m=\sum\mu_j$, or in a more traditional notation, the mixed system
of loose functional and Pfaffian equations (no rigid equations)
\beq \label{basic}
\W=0,\quad F=a; 
\eneq
defined in a certain open bounded semialgebraic subset
$$
M=M(r)\subset I_r\times B_r\times\Bbb R^{2s}
$$
{\rm(see Definition \ref{f1.1})}, such that the following holds:

$\bullet$ For any choice of the parameters $(\eps,\lb)\in B_r$ the number of
isolated $(x,y)$-solutions, denoted by $\cal B(\Cal T,\Cal S_a,{\bf f};r)$ 
of the specified basic system $(\Cal T,\Cal S_a,{\bf f};r)$ admits the estimate 
by the \asik \ system (\ref{basic}) on the manifold $M=M(r)$:
$$
\cal B(\Cal T,\Cal S_a,{\bf f};r)\leq \Cal K_{M(r)}\{\W,F;\om\}+k;
$$

$\bullet$ The forms $\W_k$ have coefficients which are polynomial in all
their arguments, and also in coordinates of the point $\Cal S_a\in
\Bbb R^{n+2s}$; the degrees of those polynomials do not exceed
$6 \mu+1$, where $\mu$ is the order of degeneracy of the corresponding
equilibrium point;

$\bullet$ The covering function $\rho(\cdot\,;r)$ for the phase space
$M(r)$ is polynomial in all its arguments and also in $r$, of the total
degree not exceeding $14k$; 

$\bullet$ The coordinate functions of the maps $F_\bt$ are explicitly
given as polynomials of the first degree on the $0$-jet space of
functions $J^0(I_r\times B_r,\Bbb R^n)$ with coefficients $\pm1$.
\ethm
The proof of this theorem is completely constructive and given in 
{\cite{IY3}}. We only point out degree estimates which 
are not given in {\cite{IY3}}.

\newpage
\begin{center}
Table 2. Separating solutions for Pfaffian systems associated with 
unfolding of elementary equilibrium points.
\end{center}
{\small
\begin{tabular}{|c |c |c |c |}
\hline
 Type
 & 
 Submanifold $\gm$
 &
 Domain  $M_r$, Covering function $\rho$
 &
 Pfaffian system $\O=0$\\ \hline
 &  &$0\ <\ x,\ y\ <\ r,$ & \\

 $S_0$ & $y=\Delta(x,\lb)$ & $\lb \in L_r,$ & $x\ dy\ - \lb y\ dx\ =0$ \\

 & & $\rho=xy(r-x)(r-y)\tilde \rho$ &  \\
 \hline 
 & & $0\ <\ x,\ y,\ z,\ w\ <\ r,$ & $x\ dz\ - \ m\ z\ dx=0,\ (1)$\\

 & $y=\Delta(x,\lb)$ & $\lb \in L_r,$ & $y\ dw\ - \ n\ w\ dy=0,\ (2)$ \\
 
 $S_\mu$& $z=x^m$ & $P_\mu(z,\lb) \neq 0,$& $m\,P_\mu(w,\lb)\,\times$ \\

 & $w=y^n$ & $\rho=xyzw(r-x)(r-y)\times$ & $y\,P_\mu(z,\lb)\,dx\,-$ \\

 & & $(r-z)(r-w)P_\mu^2(z,\lb)\tilde \rho$ & $(mP_\mu(w,\lb)+n)\,\times$   \\

 & & & $x\,P_\mu^2(z,\lb)\,dy=0\ (3)$
 \\ \hline 
 & & $|x|,|y|<r,\ x\ne0,$ & \\
 
 $D^{c}_\mu$& $y=\Delta(x,\lb)$ & $\lb\in L_r,$ & $x\ (x\ dy\ - y\ dx)=0$\\
 
 & & $\rho=(r^2-x^2)(r^2-y^2)x^2 \tilde \rho$ & \\
 \hline 
 & & $0<y<r,\ |x|<r,$ & \\
 & & $\lb \in L_r$ & \\ 
 $D^h_mu$ & $y=\Delta(x,\lb)$ & $Q_\mu(\cdot,\lb)|_{[x,1]}>0,$ &
 $Q_\mu(x,\lb)\ dy\ -y\ dx\ =0$ \\
 & & $\rho=y(r-y)(r^2-x^2)\times$ & \\
 & & $Q_\mu(x,\lb)\tilde \rho$ & \\ \hline
\end{tabular}}

$$
\ 
$$
 
{\it Notes to the Table}
Here we use the same notation as in Table 1 (and in fact this Table
continues Table 1). In particular, $n:m$ is the hyperbolicity
ratio in the resonant saddle case $S_\mu$.

In the third column of the Table the symbol $L_r$ stands for a small
$r$-cube in the $(\mu+1)$-dimensional space of the parameters 
$\lb$,
centered at the localization point $\bf c=(0,\dots,0,c)\in\Bbb
R^{\mu+1}$, corresponding to the unperturbed system:
$$
L_r=\{\lb\in\Bbb R^{\mu+1}: |\lb_i|<r,\ i=0,\dots,\mu-1,\
|\lb_\mu-c|<r \}
$$
Everywhere in the Table the function $\tilde \rho=\tilde \rho(\lb)$ is the
covering function for the set $L_r$, defined as
$$
\tilde \rho(\lb)=(r^2-\lb_1^2)\cdots(r^2-\lb_{\mu-1}^2)\cdot(r^2-(c-\lb_\mu)^2).
$$
This is a polynomial of degree $2\mu$ in all variables $\lb,r,c$.
Recall that $\deg P_mu=2\mu$ and $\deg Q_\mu=2\mu+1$ (see Table 1).
Each covering function $\rho$ is therefore a polynomial (explicitly
written in the Table). Thus, we obtain the following degree estimates:

Type $S_0$: $\deg \W=1$ and $\deg \rho =4\mu+4$.

Type $S_\mu$, $\mu>0$: $\deg \W \leq 6\mu+1$ and $\deg \rho =6\mu+8$.

Type $D^c_\mu$: $\deg \W=2$ and $\deg \rho=2\mu+6$.

Type $D^h_\mu$: $\deg \W=2\mu+1$ and $\deg \rho=4\mu+5$.

Along with the estimate  $\sum \mu_j \leq k$ (the sum of codimension)
this gives the estimates $\deg \W \leq 6\mu+1$ and $\deg \rho \leq 14k$.   

\subsection{Principal functional--Pfaffian system}

We proceed with writing down the principal functional--Pfaffian
system explicitly. Slightly abusing notation, we add the subscript $j$
for objects related to the $j$th singularity, while letters
without this subscript refer to objects related to the entire
polycycle. In this notation we omit the reference to the characteristic
size, still keeping in mind that all formulas are explicitly
polynomial in $r$.

{\bf Notations}
Denote by $M_j$ the domain from Table 2,
associated with the $j$-th singular standard map, let $\gamma_j$ be the
corresponding manifold (separating solution) and by $\W_j$ the tuple of
Pfaffian forms on it: if the singularity is of the type $D_\mu$ or
$S_0$, then $\W_j$ consists of only one form
$\om_j=A_j\,dx_j+B_j\,dy_j$, while in the case $S_\mu$, $\mu>0$, there
are three forms, of which we denote the third one by
$A_j\,dx_j+B_j\,dy_j$, (see Table 2). The covering function for $M_j$
is denoted by $\rho_j:M_j\to\Bbb R^1_+$.

{\it {Construction of the principal system}} \ \
The phase space for the principal functional--Pfaffian system is the
Cartesian product of phase spaces corresponding to all the vertices of
the polycycle and the $r$-cube in the $\eps$-space:
\beq
\begin{aligned}\label{D}
M=M(r)=M_1\times&\cdots\times M_n\times\tilde B_r, \\
M_j=M_{j,r}{\text{ are taken from the}}& {\text{ second column of Table 2}},\\
\tilde B_r=\{|\eps_i|<r,&\ i=1,\dots,k\}.
\end{aligned}
\eneq

Dimension of the phase space is equal to $2n+2s+k+m$, where:

$\bullet$ $k$ is the number of the parameters $\eps$
(the principal integer index);

$\bullet$ $n\leq k$ is the number of vertices;

$\bullet$ $s\leq n$ is the number of
resonant saddles on the polycycle (each such a vertex contributes two
additional variables $z_j,w_j$ into the list of independent variables);

$\bullet$ $m=n+\sum\mu_j\leq n+k\leq 2k$ is the number of additional 
free parameters $\lb=(\lb^1,\dots,\lb^n)$, $\lb^j\in \Bbb R^{\mu_j+1}$.

The covering function for such a space is the product
\beq \label{C}
\rho=\rho_1 \cdots \rho_n \cdot \rho_\eps : M \to \Bbb R^1_+,
\eneq
where the last factor is the covering function for $\tilde B_r$. From
Table 2 it is clear that $\rho$ is a polynomial of degree at most 
$\sum_j(6\mu_j+8)\leq 14k$ in both phase variables and the characteristic 
size $r$.

Each form on $M_j$ can be pulled back on $M$, yielding the form
which is independent of all the coordinates except for those related
to the $j$th vertex. Denote by $\W$ the union of the tuples
$\W^{(j)}$: thus $\W$ is itself the tuple of 1-forms on $M$, containing
$n+2s$ of them:
\beq 
\begin{aligned} \label{P}
\W=(\W^{(1)},\dots,\W^{(n)})=(\W_1,\dots,\W_{n+2s}),\\
\W^{(j)}=
\begin{cases}
\left\{\om_j\right\}\ \ \   {\text{if $j$ is not a resonant saddle,}}\\
\left\{\om_{j1},\om_{j2},\om_j\right\}\ \  {\text {otherwise,}}
\end{cases}
\\
{\text {where}}\ \ \  \om_{j1}=m_jx_j\,dz_j-z_j\,dx_j,\qquad &
\om_{j2}=n_jy_j\,dw_j-w_j\,dy_j,\\
\om_j=A_j\,dx_j+&B_j\,dy_j.
\end{aligned}
\eneq

Each $\gamma_j$ is a separating solution to the Pfaffian equation or
system of equations $\W^{(j)}=0$ on $M_j$, therefore the Cartesian
product
$$
\G=\gamma_1\times\cdots\times\gamma_n\times \tilde B_r
$$
is the separating solution to the Pfaffian system $\W=0$ on $M$.
Indeed, one may consider the chain of submanifolds
$$
\G_i=\gamma_1\times\cdots\times\gamma_i
\times M_{i+1}\times\cdots\times M_n\times \tilde B_r
$$
This chain possesses all the properties required by the definition of a
separating solution, see section \ref{asik}: there are no singular points 
of Pfaffian forms on all the manifolds from this chain, and the topological
condition of $\G_{i+1}$ being the boundary of a domain in $\G_i$ is
trivially satisfied, because each $\gamma_{i+1}$ is the boundary of the
corresponding subdomain in $M_j$. Thus the Pfaffian part of the
principal system is constructed.

In this Pfaffian part we have the following information about the
polynomials (recall that $\cal S_a$ stands for the algebraic part of
the specification for the basic system, which is identified by (\ref{3.0})
with a tuple of real variables):
\beq 
\begin{aligned} \label{E'}
A_j, B_j\in\Bbb Z[x,y,\lb,\Cal S_a],\qquad \deg A_j,\deg B_j\leq 6\mu+1,
\\ \rho\in\Bbb Z[x,y,\lb,\eps,r],\qquad \deg \rho \leq 14k.
\end{aligned}
\eneq

Now we proceed with description of the functional part of the
principal system. It is given by the map
\beq 
\begin{aligned} \label{F}
F=(F_1,\dots,F_{n+k+m}):M\to\Bbb R^{n+k+m},\\
F_j=
\begin{cases}
x_{j+1}-f_j(y_j,\eps),&\quad j=1,\dots,n\mod(n),\\
\eps_{j-n},&\quad j=n+1,\dots,n+k,\\
\lambda_{j-n-k},&\quad j=n+k+1,\dots,n+k+m.
\end{cases}
\end{aligned}
\eneq
The dimension of a generic fiber $F^{-1}(\cdot)$ is equal to the
codimension of separating solutions of the Pfaffian system. An
essential feature of the above map is the following one: {\it the
coordinate functions of the map $F$ are polynomial combinations of the
coordinates on the source space and generic functions $f_j$\/}:

\beq \label{E''}
F_j\in\Bbb Z[x,\eps,\lb,\bf f],\qquad \deg F_j=1,
\eneq
and all coefficients of those polynomials are $\pm1$. A more invariant
way of formulating the same property is to say that $F$ is a
polynomial map defined on the space of $0$-jets of vector-functions
\beq
\bf f:M\to\Bbb R^n,
\eneq
and this phrase makes sense since $M$ is a subset of a Euclidean space.

\bdef \label{d3.0}
The functional--Pfaffian system with the Pfaffian
equations (\ref{P}), the functional equations (\ref{F}), defined on the
domain (\ref{D}) considered with the covering function (\ref{C}), will be
called the {\it principal functional--Pfaffian system\/}. The
information provided by the estimates (\ref{E'}), (\ref{E''}) allows us to
say that the principal system is effectively described.
\endef

Later on we will refer to the principal system as simply the system
(\ref{basic}).

\subsection{ Reduction to singularity theory}

The system (\ref{basic}), whose \asik\ number majorizes the number of
solutions to the basic system (\ref{1}), satisfies the conditions of
Theorem \ref{mkh}.  The conclusion of the latter claims that the number $\Cal
K\{\W,F;\om\}$ is in turn majorized by the combination of \asik\ numbers for
some $2n+2s+k+m+1$ entirely rigid systems (recall that
$n+2s$ is the number of Pfaffian equations and $n+k+m$ is the number of 
loose functional equations in the principal system, which should
be eliminated). The properties of the principal system, listed in the
formulation of Theorem {\ref{b->f-pf}, yield a complete description of the
resulting systems as {\it chain maps\/} (the definition is given below).

In what follows we treat the original variables $x_j,y_j$, the
auxiliary variables $z_{j_\alpha},w_{j_\alpha}$ and the parameters
$\eps,\lb$ in almost the similar way, as it is suggested by the functional
equations (\ref{F}) of the principal system (\ref{basic}). The algebraic part
$\Cal S_a$ of the specification, however, plays a different role: the
coordinates of the localization points $\bf c_j$ and the integers
$n_{j_\alpha}, m_{j_\alpha}$ determining the hyperbolicity ratios of
resonant saddles, would determine the point in the new phase space,
around which the resulting chain maps will be considered.
Recall that in $\S1$ we introduced the vectors $\bf c_j$ and $\bf c$
as
\beq \nonumber
\begin{aligned}
{\bf c}_j=(0,\dots,0,c_j)\in\Bbb R^{\mu_j+1},\qquad c_j\in \Bbb R^1, \\
{\bf c}=({\bf c}_1,\dots,{\bf c}_n)\in\Bbb R^{m},\qquad m=n+\sum\mu_j.
\end{aligned}
\eneq
For our purposes it would be convenient to consider all (new) variables
as taking values around the origin in the corresponding phase space.
For this sake we make a parallel translation in the $\lb$-space, which
would take the origin into the point $\bf c$. Clearly, this
translation does not affect the algebraic structure of the principal
system (\ref{basic}), though changes the appearance of the equations.

The characteristic size $r$ retains its original meaning.

{\bf {Notations}}\ 
According to what has been said, we introduce the following notations:
\beq \nonumber
\begin{aligned}
{\bf x}&=(x,y,z,w,\eps,\lb-{\bf c}) \in \Bbb R^{2n+2s+k+m},
\\
{\bf f}&=(f_1,\dots,f_n),\qquad f_j=f_j(y_j,\eps)\iff 
{\bf f}={\bf f}({\bf x}),
\end{aligned}
\eneq
where $\bf f$ is now considered as a vector-function of the argument
$\bf x$, though each coordinate function $f_j$ of the vector $\bf
f$ depends in fact only on some of the coordinates of the vector $\bf
x$. By $D^\ell \bf f$ we denote the collection of all partial
derivatives of functions $f_j$ of the order $\ell$.

We will also use the same notation $M(r)$ for the domain of the
principal system, though in fact it would become a subset of the unit
cube $\|{\bf x}\|<r$ centered at the origin in the $\bf x$-space.

Now we can formulate the properties of the systems of equations which
appear after elimination of Pfaffian equations from the principal
system (\ref{basic}) as this was described in \S\ref{applications}.
Let ${\bf m}=2n+2s+k+m$

\bthm \label{pf->chain} \ Let ${\bf m}=2n+2s+k+m$. 
For any fixed combinatorial type $\Cal T$ of the principal
functional-Pfaffian system (\ref{basic}), any choice of the algebraic
part $\cal S_a$ of the specification and sufficiently fast decaying to zero 
sequence of numbers $\eps_1, \dots, \eps_{\bf m}$, the number of 
nondegenerate solutions to the principal system in the domain $M(r)$ for any 
choice of the characteristic size $r>0$ does not exceed the sum of the \asik\
numbers for $\bf m+1$ entirely rigid system of equations in the same domain.  
Each of these systems has the form
\beq
\begin{aligned} \label{chainmap}
{\bf P}\big({\bf x},{\bf f}({\bf x}),
D^1{\bf f}({\bf x}),\dots,D^{\bf m}{\bf f}({\bf x});\Cal S_a,r\big)
=(\eps_1, \dots, \eps_{\bf m}), \ 
{\bf x}\in M(r)\subseteq\Bbb R^{\bf m},
\end{aligned}
\eneq
where

$\bullet$ ${\bf m}\leq 7k$ is the total number of variables 
( the dimension of the phase space);

$\bullet$ ${\bf P}$ is a vector polynomial, ${\bf
P}=(P_1,\dots,P_{\bf m})$, $P_i\in\Bbb Z[{\bf x},\dots;\Cal S_a,r]$; 
the degrees  of  each polynomial $P_i$ is bounded by $14k2^i$\ $i=1,\dots, \bf m$;

$\bullet$ the domain $M(r)$ belongs to the $r$-cube of the space
$\Bbb R^{\bf m}$, centered at the origin.
\ethm

\subsection{ Chain maps and related finiteness theorems}

Now we proceed with a more invariant description of the geometric
object corresponding to the system of equations (\ref{chainmap}).

\bdef
Let $\Bbb R^{\bf m}$ be a Euclidean space with a fixed coordinate
system ${\bf x}=(X_1,\dots,X_{\bf m})$, and $U\subseteq \Bbb
R^{\bf m}$ a domain of the rectangular form,
$$
U=\{\alpha_i<X_i<\bt_i,\ i=1,\dots,\bf m\}.
$$
Denote by $I$ the index subset $I=\{1,\dots,\bf m\}$ enumerating
the coordinates in $\Bbb R^{\bf m}$, and let for any $j=1,\dots,\bf
n$ $I_j$ be a nonempty subset of $I$,
$$
\varnothing \ne I_j \subseteq I,\qquad j=1,\dots,\bf n.
$$
We say that a vector-valued function
$$
{\bf f}: U \mapsto\Bbb R^{\bf n},\qquad
{\bf f}=(f_1,\dots,f_{\bf n}),
$$
is a {\it Cartesian function\/} of the Cartesian type $\Cal
I=(I_1,\dots,I_{\bf n})$, if for any $j$ the $j$th component of this
function depends only on the coordinates $X_{i}$ with $i\in I_j$: in
other words,
$$
\forall i\notin I_j\quad \frac{\partial f_j}{\partial X_i}\equiv 0.
$$
\endef
For any given Cartesian type $\Cal I$ with ${\bf n=1}$ the set of all 
$C^\ell$-smooth
Cartesian functions (i\.e\. Cartesian maps with $\bf n=1$) of this
type constitutes a Banach space with the natural $C^\ell$-norm. We
denote this space by ${\bf C}^\ell_{\Cal I}$, sometimes omitting the
explicit reference to the type $\Cal I$ when the latter is clear from
context.  The space ${\bf C}^\ell_{\Cal I}$ will be referred to as the
{\it Cartesian space\/}. In the same way the {\it Cartesian
spaces of maps\/} arise. As a consequence, we may say about {\it
genericity\/} of Cartesian maps (functions) within the given
Cartesian type; the notions of openness and density of subsets are
also naturally defined.

\bdef
Let $\bf f$ be a $C^\ell$-smooth Cartesian map of a given Cartesian
type $\Cal I$, and $s\ge0$ an nonnegative integer number, $s\le\ell$. A
{\it Cartesian $s$-jet\/} of the function $\bf f$ at a point $\bf
x_0\in U$ is the equivalence class of all Cartesian functions of the
same Cartesian type, which differ from $\bf f$ by a term which is
$s$-flat at $\bf x_0$:
$$
\j^s{\bf f}({\bf x}_0)=\{{\bf g} \in {\bf C}^\ell_{\Cal I}:
|{\bf f}-{\bf g}|=o(|{\bf x}-{\bf x}_0|^s)\}.
$$
The space of all $s$-jets of functions of the given Cartesian type
$\Cal I$ at all points ${\bf x}_0\in U$ will be denoted by ${\bf
J}^s_{\Cal I}(\Bbb R^{\bf m},\Bbb R^{\bf n})$ or simply by $\J^s$,
when the environment is unambiguously defined by the context.

The map
$$
{\bf x} \mapsto \j^s {\bf f}(\bf x)
$$
is called the Cartesian $s$-jet extension of the Cartesian map $\bf
f$.
\endef

The space of Cartesian jets of any type and any finite order admits a
natural coordinate system, in which the Cartesian jet extension of a
map ${\bf f}=(f_1,\dots,f_{\bf n})$ takes the form
\begin{eqnarray*}
\begin{aligned}
{\bf x}=(X_1,\dots,X_M) \mapsto
\biggl( {\bf X},{\Cal F}({\bf X}), 
\left\{ \dfrac {\partial {F_j}}{\partial {X_i}},\ i\in I_j 
\right\},\dots, \\
\biggl\{ \begin{aligned} 
\text{all partial derivatives of functions}\ f_j \ 
\text{of all orders up to} \\
\ s\ \text{in the variables on which each }\  f_j \ 
\text{actually depends}
\end{aligned}\biggr\}\biggr)
\end{aligned}
\end{eqnarray*}
The Cartesian jet spaces possess almost all properties of the standard
jet spaces. In particular, the natural projections
\beq
\begin{CD}
{\Bbb R}^M\supseteq U@<pr_0<<{\bf J}^0\simeq
{\Bbb R}^M \times {\Bbb R}^{K}@<pr_1<<
\cdots@<pr_s<<{\bf J}^s@<pr_{s+1}<<\cdots
\end{CD}
\eneq
are well defined and endow each $\bf J^s_{\Cal I}$ with the structure
of an affine bundle over $\Bbb R^{\bf m}$. Thus it makes sense to say
about {\it polynomial functions\/} defined on Cartesian bundles.

\bdef
A {\it chain map\/} with the {\it exterior part\/} $\bf P$ and the
{\it interior part\/} $\bf f$ is a map of the form
$$
\Bbb R^{\bf m}\supseteq U\owns{\bf x}\mapsto {\bf P}(\j^s_{\Cal
I}{\bf f}({\bf x}))\in\Bbb R^{\bf m},
$$
where:

$\bullet$ $\bf f$ is a Cartesian map from a certain Cartesian space
$\bf C^\ell_{\Cal I}(\Bbb R^{\bf m},\Bbb R^{\bf n})$, and $\j
^s_{\Cal I}$ is the corresponding $s$-jet extension of $\bf f$;

$\bullet$ ${\bf P}:{\bf J}^s_{\Cal I}(\Bbb R^{\bf m},\Bbb R^{\bf
n})\to\Bbb R^{\bf m}$ is a vector polynomial (eventually depending
polynomially on some additional parameters),

$\bullet$ the composite map is between the spaces of 
the same dimension:$\dim{\bf x}=\dim{\bf P}={\bf m}$.
\endef

Having introduced the notions of Cartesian functions, maps, jets etc,
we can describe the system (3.10) as a chain map defined on a small
cube of some size $r>0$ with the exterior part $\bf P$ which is a
polynomial with integer coefficients and of a controlled complexity;
this polynomial depends on $r$ and some additional variables $\bf A$
as well, and the interior part $\bf f$ belongs to some Cartesian
space, since the functions $f_j$ depend only on some components of the
vector ${\bf x}=(x,y,z,w,\eps,\lb-{\bf c})$ (recall that all nonzero
coordinates of the vector $\bf c$ are already included among the
variables $\cal S_a$). Thus our problem of estimating cyclicity of a
polycycle takes the following form: describe {\bf the} Cartesian maps 
$\bf f$ for which the chain map admits an upper estimate for the 
number of preimages of regular values.

Consider chain maps of the form
\beq \nonumber
{\bf x} \mapsto {\bf G}_r(\bf x)={\bf P}(\j^s_{\Cal I}
{\bf f}({\bf x}),r)=(P_1,\dots ,P_{\bf m})(\j^s_{\Cal I}
{\bf f}({\bf x}),r), \qquad {\bf x}\in U\subset \Bbb R^{\bf m},\ r>0,
\eneq
depending polynomially on an additional variable $r$, so that
\beq \label{poly}
{\bf P}:{\bf J}^s_{\Cal I}(\Bbb R^{\bf m},\Bbb R^{\bf
n})\times\Bbb R^1\to\Bbb R^{\bf m},\quad {\bf f}\in{\bf C}^\ell_{\Cal
I}(U,\Bbb R^{\bf n}).
\eneq
We assume that the polynomial $\bf P$ and the Cartesian type 
$\Cal I$ are fixed (and $U$ denotes as before a unit cube) and 
$\bf P$ is nontrivial polynomial, i.e. at some point $x \in U$ 
the linearization matrix $dP(x)$ has full rank.

Suppose that the smoothness order $\ell$ is sufficiently high,
$$
\ell>\bf m+1.
$$

\bthm \label{chainest} 
For any polynomial 
${\bf P}=(P_1,\dots ,P_{\bf m})$ as in (\ref{poly}) one may choose a
subset 
$\mathsf F_{\bf P}\subset{\bf C}^\ell_{\Cal I}(U,\Bbb R^{\bf n})$ 
in the space of Cartesian functions of the given type, which is
open and dense in this space such that for any Cartesian function 
${\bf f}\in\mathsf F_{\bf P}$ and any sufficiently quickly decaying 
sequence $a_1, \dots ,a_{\bf m}$ there exists a characteristic
size $r_0>0$ such that the number of preimages of 
$(a_1, \dots ,a_{\bf m})$ admits the following upper estimate:
\beq
\limsup_{r\to 0^+}\#\{{\bf x}\:\ {\bf x}\in U_r,\ 
{\bf G}_r({\bf x})=
(a_1, \dots ,a_{\bf m})\} \leq \prod_{i=1}^{\bf m} \deg P_i
\leq 2^{25k^2}.
\eneq
\ethm

A bit of terminology:
``Replace an $n$-th jet $j^nF$ by its linear part at a point 
$a \in \Bbb R^n$'' means ``replace the map 
$j^nF:\Bbb R^n \to J^n(\Bbb R^n,\Bbb R^n)$ by its 
linear part $L_{F,a,n}$ at the point $a$''.

By the phrase ``a map $G:M \to N$ of manifolds satisfies 
{\it{a transversality condition}}'' we  mean
that for some manifold (resp. a collection of manifolds) in the 
image $N$ the map $G$ is transversal to this manifold (resp. 
these manifolds).

The second stage consists in constructing a stratification 
of the $n$-jet space $J^n(\Bbb R^n, \Bbb R^n)$ 
(a decomposition into a disjoint union of manifolds 
described below)
such that if the $n$-jet $j^nF$ is transversal to all
manifolds of this stratification, then 
the following theorem is true:

\bthm \label{strat} Let $P=(P_1,\dots, P_n)$ be a nontrivial
polynomial defined on the space of $n$-jets 
$P:J^n(\Bbb R^n,\Bbb R^n) \to \Bbb R^n$ and let 
$F:\Bbb R^n \to \Bbb R^n$ be a $C^k$ smooth map, $k>n$.
Suppose the $n$-jet $j^nF$ satisfies a transversality condition 
depending only on $P$. Then for a sufficiently small $r$ one can 
replace in the statement of the previous theorem
the $n$-jet $j^nF$ at the point $a$
by its linear part $L_{F,a,n}$. Namely,
\beq \label{linear}
 \begin{aligned} \label{stratif}
 \# \{x \in B_r(a):\ 
\ P_1 \circ j^nF (x)=a_1, \dots , P_n \circ j^nF(x) =a_n \}=\\
 \# \{x\in B_r(a):\ \ P_1\circ L_{F,a,n}(x) =a_1,\dots , 
 P_n \circ L_{F,a,n}(x)=a_n\},
 \end{aligned}
\eneq
where $a_1, \dots, a_n$ go to zero sufficiently fast.
By Bezout's theorem the number of solutions to the 
equation in the right-hand side of (\ref{linear}) can be bounded 
by the product $\prod_{i=1}^n \deg P_i$. 
\ethm

The classical transversality theorem
{\cite{AGV}} says that for a generic map $F$ its $n$-jet $j^nF$ 
satisfies any ahead given transversality condition.

\subsection{Stratified manifolds}

Now we recall basic definitions from the theory of stratified sets.

Let $M$ be a smooth manifold, which we call the 
{\it{ambient manifold}}. 
Consider a singular subset $V \subset M$. Roughly speaking a 
stratification of $V$ is a decomposition of $V$ into 
a disjoint union of manifolds (strata) 
$\{V_\al\}_\al$ such that {\it{strata of bigger dimension are 
attached to strata of smaller dimension in a ``regular'' way}}.

``Regular'' will obtain a precise meaning in a moment, but the most 
important property is that
{\it{transversality to a smaller stratum implies transversality 
to an ``attached'' bigger stratum.}}
Now we are going to describe the standard language
of stratified manifolds and maps of
stratified manifolds. This goes back to Whitney and
Thom {\cite{W}}, {\cite{Th}}. 

Recall the Whitney Conditions (a) and (b). Condition (a)
is similar to the notion of $a_P$-stratification due to Thom
{\cite{Th}} defined in the next subsection. We shall use 
$a_P$-stratification to prove condition ({\ref{apstra}}).

Consider a triple $(V_\bt,V_\al,x)$, where $V_\bt,\ V_\al$ are 
$C^1$ manifolds, $x$ is a point in $V_\bt$ and 
$V_\bt \subseteq \bar V_\al \setminus V_\al$.

\bdef \label{whitneya}
A triple $(V_\bt,V_\al,x)$ satisfies the Whitney (a)
condition if for any sequence of points $\{x_k\} \subset V_\al$ 
converging to a point $x \in V_\bt$
the sequence of tangent planes $T_k=T_{x_k}V_\al$ converges in 
the corresponding Grassmanian manifold of $\dim V_\al$-planes in $TM$ and
$\lim T_k=\tau \supset T_xV_\bt$.
\endef

\bdef A triple $(V_\bt,V_\al,x)$ satisfies the Whitney (b)
condition if for any two sequences of points $\{x_k\} \subset V_\al$, 
$\{y_k\} \subset V_\al$ converging to a point $x \in V_\bt$ the
sequence of ``vectors'' $\frac{y_k-x_k}{|y_k-x_k|}$ converges to a vector $v \in T_xM$
which belongs to a limiting position of 
$\lim T_{x_k}V_\al=\tau$, i.e. $v \in \tau$.

Since condition (b) is local one can think of $M$ 
as Euclidean. This explains how to interpret
the vector $\frac{y_k-x_k}{|y_k-x_k|}$. 
\endef
It is easy to show that condition (b) implies condition (a).

\bdef \label{strdef}
A locally closed subset $V$ in the ambient manifold  $M$ is called 
a {\it {stratified manifold}} (set, variety) in $M$, if it is represented as a
locally finite disjoint union of smooth submanifolds $V_\al$ of $M$, called
{\it {strata}}, of different dimensions in such a way that the closure
of each stratum consists of itself and the union of some other strata
of {\it {strictly smaller dimensions}}, and Condition (b) of
Whitney is satisfied.
\endef

Any union of submanifolds satisfying condition of this definition
\beq \label{statman}
V=\cup_\alpha V_\alpha
\eneq
is called a {\it{stratification}} of $V$, and the submanifolds $V_\alpha$
are called {\it{strata}}. A set $V$ is 
{\it{stratifiable}} if there is a ``nice'' partition into strata.
By {\it{a stratified manifold}} we mean a pair $(V,\Cal V)$ 
consisting of a manifold $V$ itself and a partition $\Cal V=\{V_\alpha\}$.

\subsection{Stratified maps and $a_P$-stratification}

Now we define a smooth map of a stratified manifold $(V,\Cal V)$:
\bdef \label{stratdef}
Let $(V,\Cal V)$ be a stratified manifold in an ambient
manifold $M$, $V \subseteq M$, then a map $f:V \to N$ 
is called $C^2$-smooth if it can be extended to a $C^2$ smooth map
of the ambient manifold $M$
$F:M \to N$ whose restriction to $V$ coincides
with $f$.

A stratification $V=\cup_\alpha V_\alpha$ stratifies a smooth 
map $f:V \to \Bbb R^k$ if the restriction of $f$ to any stratum 
$V_\alpha$ has constant rank, i.e., rank $df|_{V_\alpha}(x)$ is 
independent of $x \in V_\alpha$.

A map $G:L \to M$ is called transversal to a stratified set  
$(V,\Cal V)$ if\ $G$ is transversal to each strata $V_\al \in \Cal V$.
\endef

By the Rank Theorem, if a stratification 
$(V,\Cal V),\ \Cal V=\{V_\alpha\}_{\alpha \in I}$ 
stratifies a smooth map $P$, then for each strata 
$V_\al$ the number  
$d_\alpha(P)=\dim V_\alpha- {\text{rank}}\ dP|_{V_\alpha}$
is well defined.

Assume $d_\al(P) \geq d_\bt(P)$ for each $V_\bt \subseteq 
\bar V_\al \setminus V_\al$, i.e. nonempty level sets inside
the bigger stratum $V_\al$ have dimension $d_\al(P)$
greater or equal to dimension of the level sets $d_\bt(P)$
in the smaller stratum $V_\bt$.
We require that for any sequence of points
$\{a_k\} \subset P(V_\al)$ converging to a point $a \in P(V_\bt)$,
the nonempty level sets $\{P^{-1}(a_k) \cap V_\al\}$
approach the limiting level set $\{P^{-1}(a) \cap V_\bt\}$
``regularly''.
In other words, we {\it{require that the level sets in the 
bigger stratum $V_\al$
approach the limit level set in the smaller stratum $V_\bt$ nicely.}}

\bdef \label{a-p}Let $P: M \to N$ be a $C^2$ smooth  map of manifolds, 
and let $V_\bt$ and $V_\al$ be submanifolds of $M$ such that the restrictions
$P|_{V_\bt}$ to $V_\bt$ and $P|_{V_\al}$ to $V_\al$ have constant ranks 
$R_{V_\bt}(P)$ and $R_{V_\al}(P)$, respectively.  
Let $x$ be a point in $V_\bt$.

We call the manifold $V_\al$ $a_P$-regular over $V_\bt$ with respect to 
the map $P$ at the point $x$ if for any sequence of points 
$\{x_n\} \subset V_\al$
converging to $x \in V_\bt$ the sequence of tangent planes
to the level sets $T_k=ker\ dP|_{V_\al}(x_k)$ converges 
in the corresponding Grassmanian manifold of 
$(\dim V_\al-R_{V_\al}(P))$-dimensional planes to a plane $\tau$ and 
\beq \label{apcond}
\lim ker\ dP|_{V_\al}(x_k)=\tau \supseteq ker\ dP|_{V_\bt}(x)
\eneq
\endef

\bdef \label{apstratif}
A $C^2$ smooth map $P:V \to N$ of a stratifiable
manifold $V$ to a manifold $N$ is called 
$a_P$-stratifiable if there exist a stratification 
$(V,\Cal V)$ such that the following 
conditions hold:

a) $(V,\Cal V)$ stratifies the map $P$ (see definition \ref{stratdef});

b) for all pairs $V_\bt$ and $V_\al$ from $\Cal V$ such that 
$V_\bt \subseteq \bar V_\al \setminus V_\al$ the stratum $V_\al$ is
$a_P$-regular over the stratum $V_\bt$ with respect to $P$ at point $x$
for all $x \in V_\bt$.
\endef

The original definition of $a_P$-stratification 
requires an appropriate stratification of the image also
{\cite{Ma}}, but we do not require stratification of the image for
our purposes.

\subsection{Relation between existence of $a_P$-stratification and
condition (\ref{apstrat}).}

In section \ref{heur} we showed that the key to the proof
of Theorem \ref{chainest} is condition 
(\ref{apstra}) (see Proposition \ref{lin}).
Now we are going to reduce the question whether
condition (\ref{apstra}) is satisfied to the question
whether an $a_P$-stratification of the polynomial $P$ exists.

Let $P=(P_1,P_2):\Bbb R^N \to \Bbb R^2$ be a nontrivial 
polynomial, $V=P_2^{-1}(0)$ and  $V_0=(P_1,P_2)^{-1}(0)$ be level 
sets. Assume that there exists a stratification
$(V,\Cal V)$ that stratifies the map $P|_V$ such that the zero level set $V_0$ 
can be represented as a union of strata from 
$\Cal V$, i.e., $V_0=\cup_{\alpha \in I_0} V_\alpha$.
Denote this stratification of $V_0$ by $\Cal V_0$. 
Recall that a map $F:\Bbb R^k \to \Bbb R^N$ is transversal to  a
stratification $(V_0,\Cal V_0)$ if it is transversal to each
strata $V_\al \in \Cal V_0$.
Associate to each level set $V_a, \ a \neq 0$ a natural 
decomposition $\Cal V_a=\{V_a \cap V_\alpha \}_{\alpha \in I}$.

\bprop \label{app} With the above notation  
if a stratum $V_\al \in \Cal V \setminus \Cal V_0$ is
$a_P$-regular over a stratum $V_\bt \in \Cal V_0$ with respect to 
the polynomial $P$, then any $C^2$ smooth map 
$F: \Bbb R^n \to \Bbb R^2$ transversal to $(V_0,\Cal V_0)$
 is also transversal to $V_a \cap V_\alpha$ for any small 
$a$. This is equivalent to condition (\ref{apstra}). 
\enprop 

{\it{Proof}}\ \ 
Pick a point $x$ in $V_\bt \subset V_0$ and a point $y \in V_\al$. 
Notice that $ker\ dP|_{V_\bt}(x)$ is the tangent 
plane to the level set $\{P^{-1}(P(x)) \cap V_\bt\}$ 
at the point $x$ and $ker\ dP|_{V_\al}(y)$ is the tangent plane
to the level set $\{P^{-1}(P(y)) \cap V_\al\}$.

By condition (\ref{apcond}) if a map $F:X \to \Bbb R^N$ is transversal to 
$ker\ dP|_{V_\bt}(x)$ at a point $x$, then $F$ is transversal to
$ker\ dP|_{V_\alpha}(y)$ for any $y \in V_\alpha$ near $x$.

Therefore, the condition ``$F$ is transversal  to $V_\bt$ at a point 
$x$''
implies the condition ``$F$ is transversal to $V_\alpha \cap V_a$ 
for any small $a$''.
This completes the proof. 

\subsection{Existence of $a_P$-stratification
for polynomial maps}

The existence of $a_P$-stratifications is not a trivial question.
There are some obvious obstacles. For example,
let $V \subset \Bbb R^n$ be an algebraic variety and 
let $P:\Bbb R^n \to \Bbb R^k$ be a polynomial map.
Assume that $(V,\Cal V)$ stratifies $P$.
If we have two strata $V_\alpha$ and $V_\bt$ so that $V_\al$ lies ``over''
$V_\bt \subseteq \bar V_\alpha \setminus V_\alpha$, then 
condition (\ref{apcond}) can't be satisfied if dimension 
of the level sets $d_\al(P)$ in the upper stratum $V_\alpha$ is strictly
less than that of $d_\bt(P)$ in the lower stratum $V_\bt$, 
i.e., $\dim ker\ dP|_{V_\alpha}(y)< \dim ker\ dP|_{V_\bt}$.
In this case a plane $ker\ dP|_{V_\bt}(x)$ of the
lower stratum $V_\bt$ should belong to a plane 
$\tau$ of smaller dimension (see condition (\ref{apcond})), 
which is impossible. Thom constructed the first example when 
this happens {\cite{GWPL}}.

{\it{Thom's example}}

Consider the vector-polynomial $P$ in the form $P:(x,y) \to (x,xy)$.
The line $\{x=0\}$ is the line of critical points of $P$.
Outside of the line $\{x=0\}$ $P$ is a diffeomorphism. Therefore, 
the preimage of any point $a \neq 0$ $P^{-1}(a)$ is $0$-dimensional.
On the other hand, the preimage of $0$ is the line $\{x=0\}$.

\bdef Let us call an algebraic set $V$ rank compatible with 
a polynomial $P$ if there exists
a stratification $(V,\Cal V)$ which stratifies $P$ and 
for any pair $V_\alpha$ and $V_\bt$ from $\Cal V$ such that 
$V_\bt \subseteq 
\bar V_\al \setminus V_\al$ dimensions of the levels 
$d_\bt(P)$ in the lower stratum $V_\bt$ do not
exceed dimensions of the level sets $d_\alpha(P)$ in the upper 
stratum $V_\alpha$.
\endef 
It turns out that even if an algebraic set $V$ is rank 
compatible with a polynomial $P$, then $a_P$-stratification 
still does not always exist. Let us present an example
with this property. The example below belongs to M.Grinberg.
It seems that the existence of a 
counterexample was known before, but we did not find an 
appropriate reference.

\subsection{Nonexistence of $a_P$-stratification}
Let $V=\{(x,y,z,t)\in \R^4:\ x^2=t^2y+z\}$ be the
three dimensional algebraic variety and
$P:V\to \R^2$ be the natural projection to the last 
two coordinates, i.e. $P:(x,y,z,t)\to (z,t)$.

\bprop \label{Grinb} With the above notations the set $V$ is rank
compatible with the polynomial map $P$ and
does not have $a_P$-stratification.
\enprop

{\it{Proof}} \ \ Consider a rank stratification of $V$.
Such a stratification consists of three stratum:
$V_1=\{x=t=z=0\},\ V_2=\{t=0,\ x^2=z, \ x \neq 0\},$ and 
$V_3=\{t \neq 0\}.$  On each stratum  rank$P|_{V_i}=i-1$.
Level sets $P^{-1}(t,z)$---parabolas for $t \neq 0$ and
lines for $t=0$. 

Show that for each point ${\bf a}=(0,a,0,0) \in V_1$
there exists a family of level sets such that at the point
$\bf a$ the property $a_P$-regularity of $V_3$ over $V_1$ fails.

Consider the preimage of the curve $\{z=-at^2\} \subset \Bbb R^2$.
This is an algebraic variety of the form $W_a=\{x^2=t^2(y-a)\}$.
One can see that $W_a$ is the Whitney umbrella.©The level
$x^2=t_0^2(y-a)$ is the parabola. As $t_0 \to 0$ this
parabola tends to semiline ${x=t=z=0, y \geq a}$.
At the point $\bf a \in V_1$ the property $a_P$-regularity of 
$V_3$ over $V_1$ clearly fails. This completes 
the proof of the Proposition.

Let us mention a positive result on existence of 
$a_P$-stratification.

\bthm \cite{Hir1} \label{Hir}
If $V \subset \Bbb R^n$ is a semialgebraic variety and
$P:\Bbb R^n \to \Bbb R$ is a polynomial function, then 
there exists an $a_P$-stratification of $(V,\Cal V)$
with respect to $P$.
\ethm

\section{Existence of $a_P$-stratification.}
\label{apstrat}
In this section we prove existence of $a_P$-stratification
in the special case we are interested in. As the Example
\ref{Grinb}  shows, the existence of a $a_P$-stratification
is a nontrivial question. In general, it does not exist.
Unfortunately, the existence of a $a_P$-stratification in our 
case does not follow from the classical results, so we need 
to prove it. 

Let $\Bbb R^N$ and $\Bbb R^k$ be Eucledian spaces with the 
fixed coordinate systems $x=(x_1,\dots, x_N) \in \Bbb R^N$ and 
$a=(a_1,\dots, a_k) \in \Bbb R^k$ with $N \geq k$ and a 
non-trivial vector-polynomial $P:\Bbb R^N  \in \Bbb R^k$.
Recall that $P$ is a nontrivial if it has a point $x \in \R^N$,
where rank $dP(x)=k$. In what follows we call {\it{vector-polynomial}} 
by {\it{polynomial}} for brevity.  

\bdef
Let ${\bf m}=(1,m_2, \dots, m_{k}) \in \Z^{k}_+$ 
and $\dt > 0$. We call the $(\bf m, \dt)$-cone $K_{{\bf{m}},\dt}$
the following set of points
\beq
\begin{aligned}
K_{{\bf m}, \dt} = \{a=(a_1, \dots, a_k) \in 
\R^k:\  0 < a_1 < \dt,\\
0 < |a_{j+1}| < |a_1 \dots a_j|^{m_{j+1}}\ 
\text{for}\ j=1, \dots, k-1\}.
\end{aligned}
\eneq
Let $\m'=(1,m'_2, \dots, m'_{k})\in \Z^{k}_+$. Define 
$\m'\succ \m$ if $\m' \neq \m$ and $m'_i\geq m_i$ for all
$2 \leq i \leq k$. We call the $(\bf m', \dt')$-cone $K_{{\bf{m'}},\dt}$  
a refinement of the $(\bf m, \dt)$-cone $K_{{\bf{m}},\dt}$ if
$\m'\succ \m$ and $\dt\geq \dt'$.
\endef

Define the following sets
\beq \label{limset}
V_{{\bf m}, \dt, P} = 
\textup{closure} \{ P^{-1}(K_{\bf m, \dt})\},\  
V_{0, \bf m, P} = \cap_{\dt > 0} V_{\bf m, \dt, P}
\eneq
Then one has 
\bthm \label{existence} For any nontrivial polynomial $P$
there exist an integer vector ${\bf{m}} \in \Bbb Z_+^n$ and 
positive $\dt$
such that the following conditions hold

a) the set  $V_0=V_{0,{\bf{m}},P}$ (see (\ref{limset})) is 
semialgebraic. 

b) the set  $V_{{\bf{m}},\dt,P}$ consists of regular 
points of $P$, i.e. if $b \in V_{{\bf{m}},\dt,P}$, then 
the level set $P^{-1}(b)$ is a manifold of codimension $n$.

c) there exists a stratification of $V_0$ by semialgebraic 
strata $(V_0,\Cal V_0)$ satisfying the property: 
$V_{{\bf{m}},\dt,P}$ is $a_P$-regular over 
any strata $V_\alpha \in \Cal V_0$ with respect to $P$.
\ethm

In order to prove Theorem \ref{existence} we reformulate 
it in a convenient for us language. Let $a \in \R^k$. 
Denote by $L_a=P^{-1}(a)$ the level set of $P$. 
Recall that $a \in \R^k$ is called a regular value if 
for any $x \in L_a$ the rank of linearity of $P$ is maximal, 
i.e. rank$dP(x)=k.$
 
\bdef
Let $a,b \in \R^k$ be values of $P:\ \R^N \in \R^k,\  
B^N \subset \R^N$ be the unit ball centered at the origin, and 
\beq
\begin{aligned}
d_0^0: B^N \times \R^k \in \R,&\  d_0^0(x,a) = 
\inf_{y \in L_a \cap B^N} \|x-y\|^2
\\
d_0: \R^k \times \R^k \in \R,&\ d_0(a,b) = 
\sup_{x \in L_b \cap B^N} d_0^0(x,a).
\end{aligned}
\eneq
Then the $C^0$-distance between level sets $L_a \cap B^N$ and 
$L_b \cap B^N$
$$
D_P^0(a,b)=dist_C^0(L_a \cap B^N, L_b \cap B^N)
=\frac 12 \left(d_0(a,b)+d_0(b,a)\right).
$$
\endef

For any $1\leq m \leq N$ denote by $G^{m,N}$ the set of 
$m$-dimensional planes in the $N$-dimensional Euclidean space.
$G^{m,N}$ is so-called the grassmanian manifold. Below
we introduce convenient for as distance in the grassmanian 
manifold $G^{m,N}$. Now we define $C^1$-distance between regular 
level sets $L_a$ and $L_b$ in an appropriate for us way. Write $P$ 
using coordinate functions $P=(P_1, \dots, P_k):\R^N \to \R^k$. 
If $x \in \R^N$ is a regular point of $P$, then gradients 
$\nabla P_1(x), \dots, \nabla P_k(x)$ are linearly independent 
and span the space which is the orthogonal 
complement to the tangent space to the level set $P^{-1}(P(x))$. 
Define a Gramm-Schmidt orthogonalization operator:

\bdef \label{gramm}
Let $v_1, \dots, v_k \in \R^k$ be linear independent 
vectors. Define the Gramm-Schmidt linear operator by
\beq
\begin{aligned}
*:(v_1,\dots, v_k)=(v_1^*,\dots, v_k^*),\  \ {\textup{where}}\\ 
v_1^*=v_1,\ v_2^*=v_2 - \frac{(v_2, v_1)}{(v_1,v_1)}v_1, \dots, 
v_k^*=v_k-\sum_{j<k}\frac{(v_k, v_j)}{(v_j,v_j)}v_j.
\end{aligned}
\eneq
\endef
{\bf Remarks}\ {\it{0. The Gramm-Schmidt linear operator $*$
has nothing in common with the asterisk operator used for the
Khovanski reduction procedure in section \ref{asik}.

1. Vectors $\{v_1,\dots, v_k\}$ and 
$\{v_1^*,\dots, v_k^*\}$ span the same $k$-dimensional 
space denoted by $L$; 

2. Vectors $\{v_1^*,\dots, v_k^*\}$ form an orthonormal
basis in the plane $L$;

3. Let $\{L_t\}_{t \in (0,1]}$ be a family of
$k$-dimensional planes in $\R^N$ spanned by a family of 
vectors $\{v_1(t),\dots, v_k(t)\}_{t \in (0,1]}$ depending 
continuously on $t$. Consider 
$\{v^*_1(t),\dots, v^*_k(t)\}_{t \in (0,1]}=
*(v_1(t),\dots, v_k(t))$ as the family of orthonormal basis
in $\{L_t\}_{t \in (0,1]}$. Then sufficient condition
that $L_t \to L$ in the grassmanian manifold $G^{k,N}$
is existence of an orthonormal basis
$\{v_1^*,\dots, v_k^*\}$ of $L$ such that 
\beq
\frac{\left(v^*_j(t),v^*_j\right)^2}
{(v^*_j(t),v^*_j(t))(v^*_j,v^*_j)} \to 0 \ \ 
\textup{as} \ t \to 0.
\eneq}}

Define the Gramm-Schmidt operator for the polynomial map $P$  
\beq \label{schm}
(*(dP)_1(x),\dots,*(dP)_k(x))=
*(\nabla P_1(x),\dots ,\nabla P_k(x)).
\eneq
Each vector $*(dP)_j(x)$ is given by the rational function 
in $x$. 

Let $\Sigma \subset \R^N$ be the set of critical points of $P$. 
To measure $C^1$-distance between two regular level sets 
we introduce the following function: Let
$x,y \notin \Sigma$. Then 
\beq
\begin{aligned} \label{scalar}
R_{P}(x,y)=\sum_{j=1}^k
\left(1-
\frac{\left(*(dP)_j(x),*(dP)_j(y)\right)^2}
{(*(dP)_j(x),*(dP)_j(x))(*(dP)_j(y),*(dP)_j(y))}\right)\\
Q_{P}(x,y)=\|x-y\|^2+R_P(x,y).\quad \quad \quad 
\end{aligned}
\eneq
\bdef
Let $a,b \in \R^k$ be regular values of $P:\R^N \to \R^k$, and
$L_a=P^{-1}(a)$ and $L_b=P^{-1}(b)$ be regular level sets.
\beq \nonumber
\begin{aligned}\label{dist1}
d_{1,P}^0:B^N \times \R^k\setminus \Sigma \to \R,\quad &
d^0_{1,P}(x,a)=\inf_{y \in L_a \cap B^N} 
Q_{P}(x,y),\\
d_{1,P}:\left(\R^k\setminus \Sigma \right) \times 
\left(\R^k\setminus \Sigma\right) \to \R,  &\quad
d_{1,P}(a,b)=\sup_{x \in L_b \cap B^N} d^0_{1,P}(x,a).
\end{aligned}
\eneq
Then the $C^1$-pseudodistance between regular level sets 
$L_a \cap B^N$ and $L_b \cap B^N$ is defined by
\beq \label{dist2}
D^1_P(a,b)=dist_{C^1}\left(L_a \cap B^N, L_b \cap B^N\right)=
\frac 12\left(d_{1,P}(a,b)+d_{1,P}(b,a)\right).
\eneq
\endef
\brm We call the function $D^1_P(a,b)$ $C^1$-pseudodistance, not
$C^1$-distance, because it does not satisfy the triangle inequality.
However, it satisfies the following triangle-like inequality
\beq \label{pstriang}
2(D^1_P(a,b)+D^1_P(b,c))>D^1_P(a,b).
\eneq
The reason we define $C^1$-pseudodistance $D^1_P(a,b)$ in such a way
is because  the function $D^1_P(a,b)$ is algebraic (see
Lemma \ref{polynom} below).

The inequality (\ref{pstriang}) can be proven as follows.
Let $v,w \in \R^N$ be vectors. Denote by $\angle (v,w)$
the angle between $v$ and $w$. Direct calculation shows that
$$
R_P(x,y)=\sum_{j=1}^k \sin^2(\angle(*(dP)_j(x),*(dP)_j(y))).
$$ 
It is easy to check that $2 (\sin^2 \al +\sin^2 \bt)\geq
\sin^2 (\al +\bt)$ which is sufficient for the proof of the
inequality (\ref{pstriang}) . 
\erm

Now we can reformulate Theorem \ref{existence} in the following 
way
\bthm \label{existence1}
For any nontrivial polynomial $P$
there exist an integer vector 
${\bf{m}}=(1,m_2,\dots,m_k) \in \Bbb Z_+^n$ and 
positive $\dt$ such that the following conditions hold

a) for any two values with the same first coordinate 
$a=(t,a_2,\dots,a_k)$ and $b=(t,b_2,\dots,b_k)$ in $K_{\m,\dt}$ 
$$
D_P^0(a,b)<t;
$$

b) the same as in Theorem \ref{existence};

c) for any two values with the same first coordinate 
$a=(t,a_2,\dots,a_k)$ and $b=(t,b_2,\dots,b_k)$ in $K_{\m,\dt}$ 
$$
D_P^1(a,b)<t;
$$
\ethm

Let us show that parts a) and c) imply parts a) and c) of Theorem 
\ref{existence} respectively.

Proof a) $\Longrightarrow$ from a) of Theorem \ref {existence}.
Consider an algebraic curve of the form 
$\gm (t)=(t, t^{m_2+1}, t^{(m_2+2)(m_3+1)}, \dots, 
t^{(m_2+2)\dots (m_{k}+1)})$. One can check that  
$\gm (t) \in \ K_{\m, \dt}$ for any $t \in (0, \dt)$. 
Denote by $V_{t,P} = P^{-1}(\gm (t))$. The set 
$\cup_{0 < t \leq \dt} V_{t,P}$ is clearly 
semialgebraic set. By the Tarski-Seidenberg theorem the 
following set 
$$
\textup{closure} \{\cup_{0 < t \leq \dt} V_{t,P}\} \setminus 
\{\cup_{0 < t \leq \dt} V_{t,P}\}
$$ 
is semialgebraic. Since for any smooth curve 
$\gm'(t)= (t,\gm_2{(t)},\dots,\gm_k{(t)})\in K_{\m, \dt},\ 
t\in (0,\dt)$  Hausdorff distance  between the 
level sets $V_{t,P}=P^{-1}(\gm(t))$ and 
$V'_{t,P}=P^{-1}(\gm'(t))$ is at most $t$, i.e.  
$D_P^0(\gm(t), \gm'(t)) < t$. It implies that 
 Hausdorff distance between any two level sets
of the form $V_{t,P}$ and $V'_{t,P}$ 
tends to $0$ as $t \to 0$. Therefore, 
$$
\textup{closure} \{\cup_{0 < t \leq \dt} V_{t,P}\} \setminus 
\{\cup_{0 < t \leq \dt} V_{t,P}\} = 
\textup{closure} \{P^{-1}(K_{\m, \dt})\} \setminus 
\{P^{-1}(K_{\m, \dt})\}. 
$$
This completes the proof of part a). 

Proof c) $\Longrightarrow$ from c) of Theorem \ref{existence}.
Let us use notations of the proof of part a). 
By theorem \ref{Hir} there is a stratification of $V_{0, \m, P}$ 
such that the semialgebraic set 
$\{\cup_{0 < t \leq \dt} V_{t,P}\}$ is a 
$a_P$-regular over $V_{0, \m, P}$. Indeed, let 
$\pi_1: \R^k \to \R$ be the natural projection onto the 
first coordinate. Then a polynomial function 
$p=\pi_1 \circ P: \R^N \to \R$ is well-defined and 
$p^{-1}(t)=P^{-1}(\gm(t))$. Application of theorem \ref{Hir}
to the map 
$$
p:\textup{closure} \{\cup_{0 < t \leq \dt} V_{t,P}\}\to \R
$$
gives existence of a required stratification.

Since $C^1$-distance between any 
two level sets of the form $V_{t,P}= P^{-1}(\gm(t))$ and 
$V'_{t,P}= P^{-1}(\gm'(t))$ is at most $t$, i.e. 
$D^1_P(\gm(t), \gm'(t)) < t$. It implies  that 
$C^1$-distance between any two level sets
of the form $V_{t,P}$ and $V'_{t,P}$ 
tends to $0$ as $t \to 0$. Therefore, 
$a_P$-regularity of $P^{-1}(K_{\m, \dt})$  over 
$V_{0, \m, P}$  follows from $a_P$-regularity of 
$\{\cup_{0 < t \leq \dt} V_{t,P}\}$  over 
$V_{0, \m, P}$. This completes the proof of part c).

Before proving Theorem \ref{existence1} let us formulate 
a basic fact from elimination theory {\cite{Mu}}.

\subsection{Elimination theory}

Let $\C^m$ denote the $m$-dimensional complex space 
$z=(z_1,\dots,z_m) \in \C^m,\ \ m \in \Z_+$.
A set $V$ in $\C^m$ is called {\it{a closed algebraic set}} 
in $\C^m$ if there is a finite set of polynomials 
$F_1,\dots F_s$  in $z_1,\dots,z_m$ such that
\beq \nonumber
V(F_1,\dots ,F_s)=\{(z_1, \dots ,z_m) \in \C^m |\ 
F_j(z_1,\dots,z_m)=0,\ 1\leq j \leq s\}.
\eneq
One can define a topology in $\C^m$, called the 
{\it{Zariski topology}}, whose closed sets are closed
algebraic sets in $\C^m$. This, indeed, defines a topology,
because the set of closed algebraic sets is closed under a
finite union and an arbitrary intersection. 
Sometimes, closed algebraic sets are also called
Zariski closed sets. 

\bdef A subset $S$ of $\C^m$ is called constructible
if it is in the Boolean algebra generated by the 
closed algebraic sets; or equivalently if $S$ is a disjoint
union $T_1 \cup \dots \cup T_k$, where $T_i$ is locally
closed, i.e. $T_i=T'_i-T''_i$, $T'_i$ --- a closed algebraic 
set and $T''_i\subset T'_i$ --- a smaller closed algebraic.
\endef

One of the main results of Elimination theory is 
the following 

\bthm ({\cite{Mu}, ch.2.2}) \label{elim} Let 
$V \subset \C^\mu \times \C^N$ be 
a constructible set and $\pi:\C^\mu \times \C^N \to \C^\mu$
be the natural projection. Then  $\pi(V)\subset \C^\mu$ is a 
constructible set. 
\ethm 

\section{Proof of Theorem \ref{existence1}}
\subsection{Existence of the 
$(\bf m, \dt)$-cone $K_{\m, \dt}$ of regular values
of $P$ (or Proof of Part a) of Theorem \ref{existence1}).}

The set of critical values $\Sigma_P \subset \R^k$ 
of a nontrivial polynomial map $P:\R^N \to \R^k$ is an 
algebraic set of positive codimension. It follows from
Sard's lemma for algebraic sets {\cite{Mu}}. 
Suppose $d: \R^k \to \R$ is a nonzero polynomial 
whose zero level set $d^{-1}(0) \supseteq \Sigma$. 
Fix coordinate systems in $\R^N$. By writing the linearization 
matrix $dP: \R^N \to \R^k$ and considering $(N-k+1)$ different 
$k\times k$ minors one can calculate $d$ explicitly.

\blm \label{cone}
For a nonzero polynomial $d:\R^k \to \R$ there exists
an integer vector $\m=(1,m_2,\dots,m_k)\in \Z^k_+$ and 
$\dt>0$ such that $d$ does not vanish on the $(\m, \dt)$-cone
$K_{\m, \dt}$.
\elm
\brm
If $\Sigma_P \subseteq d^{-1}(0)$ and  
$\Sigma_P \cap K_{\m, \dt} \neq \emptyset$, then 
there exists $x \in K_{\m, \dt}$ such that $d(x)=0$. 
This shows that part a) of Theorem \ref{existence1}
follows from this Lemma.
\erm
{\it{Proof}} \ \  
 Let us prove the statement by induction in dimension $k$. 

For $k=1$ the level set $d^{-1}(0)\subset \R$ is a finite 
collection of points and Lemma is obvious.

Without loss of generality assume $d(x_1,\dots, x_k)$ 
is not divisible by $x_k$. If $d$ is divisible by $x_k$, 
then for some $\bt \in \Z_+$  one can decompose 
$d(x_1,\dots, x_k)=x_k^\bt \hat d(x_1,\dots, x_k)$
so that $\hat d(x_1,\dots, x_{k-1},0)$
is not identically zero. If for some $\m\in \Z^k_+$
and $\dt>0$ the $(\m, \dt)$-cone $K_{\m, \dt}$
does not intersect zero locus $\hat d^{-1}(0)$,
then $K_{\m, \dt}$ does not intersect zero locus 
$d^{-1}(0)=\hat d^{-1}(0) \cup \{x_k=0\}$ too.

With the assumption of indivisibility by $x_k$ 
the following set 
$\Sigma^{k-1}_P=d^{-1}(0)\cap \{x_k=0\}\subset \R^{k-1}$
is of a positive codimension.
By inductive hypothesis there exist an integer vector
$\m_{k-1}\in \Z^{k-1}_+$ and $\dt_{k-1}>0$ such
that the $(\m_{k-1},\dt_{k-1})$-cone 
$K^{k-1}_{\m_{k-1},\dt_{k-1}} \subset \R^{k-1}$
has empty intersection with  $\Sigma^{k-1}_P$.

Let $\alpha=(\al_1,\dots ,\al_k) \in \Z^k_+$ and
$|\al|=\sum_j \al_j$. Write $d(x)=\sum_{\al \in \Z^k_+}a_\al x^\al$
Denote by $\deg d=\max\{|\al|: \al\in \Z^k_+,\ a_\al \neq 0\}$. 
Put $m_k=\deg d+1$. 

\bprop \label{step}
With the above notations there exists $\dt>0$ such that
for $\m=(\m_{k-1},m_k)\in \Z_+^k$ the $(\m,\dt)$-cone
$K_{\m, \dt}$ does not intersect zero locus 
$d^{-1}(0)$.
\enprop

{\it{Proof}} \ \ Put $d(0)=0$ otherwise the proposition
is trivial. Write 
$$
d(x_1,\dots,x_k)=\sum_{j=0}^{\deg d}
a_j(x_1,\dots,x_{k-1})x_k^j.
$$
Without loss of genericity one can assume
that $a_0(x_1,\dots,x_{k-1})$ does not vanish on
the $(\m_{k-1},\dt_{k-1})$-cone $K^{k-1}_{\m_{k-1},\dt_{k-1}}$.
If not, then apply Lemma \ref{cone} and refine 
$K^{k-1}_{\m_{k-1},\dt_{k-1}}$ to a required size.
By the definition of the $(\m,\dt)$-cone $K_{\m,\dt}$
the condition  $(x_1,\dots ,x_k) \in K_{\m,\dt}$ implies that
$(x_1,\dots ,x_{k-1}) \in K^{k-1}_{\m_{k-1},\dt_{k-1}}$
and $0<x_k<(x_1\dots x_{k-1})^{m_k}$. Put
$x_k=\lb (x_1\dots x_{k-1})^{m_k}$ with $\lb \in (0,1)$.
It is easy to check that
\beq \label{prod}
d(x_1,\dots ,x_k)=
a_0(x_1,\dots,x_{k-1})(1+p(x_1,\dots,x_{k-1},\lb)),
\eneq
where $p$ is such   a polynomial that $p(0,\lb)\equiv 0$.
Indeed, the  choice of $m_k$ is such that
$(x_1\dots x_{k-1})^{m_k}=a_0(x_1,\dots,x_{k-1})
q(x_1,\dots,x_{k-1})$ for some polynomial
$q(x_1,\dots,x_{k-1})$. Since $p(0,\lb)\equiv0$
for a sufficiently small $\dt$, any
$(x_1,\dots,x_{k-1})\in K^{k-1}_{\m_{k-1},\dt_{k-1}}$,
and any $\lb \in [0,1]$ the following inequality holds
$|p(x_1,\dots,x_{k-1},\lb)|<1/2$.
This shows that $d$ does not vanish on 
the $(\m,\dt)$-cone $K_{\m,\dt}$ and completes the proof
of the Proposition. 

As we pointed out above the Proposition implies 
Lemma \ref{cone}.

\subsection{Reduction to an optimization problem 
(or Proof of parts a) and c) of Theorem \ref{existence})}

Let $P=(P_1, \dots, P_k):\R^N \to \R^k$ be a nontrivial polynomial 
with $N \geq k$ given by its coordinate functions and an $(\m, \dt)-$cone 
$K_{(\m, \dt)} \subset \textup{Im}\ P(\R^N) \subset \R^k$ be a cone of 
regular values of $P$. Existence of such a cone is proven in the 
previous section. Recall that $\Sigma \subset \R^N$ denotes the  set of critical 
points of $P$ and $Q_P(x,y)$ is defined in (\ref{scalar}). 
The function
$Q_P(x,y)$ is a rational function symmetric with respect to permutation of 
$x$ and $y$. It is defined to measure $C^1$-distance between level sets
(see remarks after definition \ref{gramm}).
The singular set of $Q_P$ belongs to  
$(\Sigma \times \R^N) \cup (\Cal \R^N \times \Sigma)$. Recall that
$B^N = \{x:\ \sum_{i=1}^N x_i^2 \leq 1 \} \subset \R^N$. Introduce 
functions $r(x) = 1 - \sum_{i=1}^N x_i^2.$ 

Assume that the restriction of P to the boundary 
$S^N = \partial B^N = \{x:\ \|x\|=1 \}$ has only the regular values 
in the $(\m, \dt)-$cone $K_{(\m, \dt)}$. Indeed, regularity of 
$P|_{S^N}:\ S^N \to \R^k$ is equivalent to regularity of the polynomial
map $(P,r):\ \R^N \to \R^k \times \R$ given by 
$(P,r)(x)=(P(x), \sum_{i=1}^N x_i^2)$.
Existence of an $(\m, \dt)-$cone of regular values of the map $(P,r)$
follows from Lemma \ref{cone}.

\blm \label{polynom}
With the notations above let 
$a_{\tau_1},\ a_{\tau_2} \in K_{(\m, \dt)}$ be two points with the
same first $(k-1)$ coordinates, i.e. $a_{\tau_1}=(a^{k-1},\tau_1)$
and $a_{\tau_2}=(a^{k-1},\tau_2)$. Then there exists a polynomial
$R(a^{k-1},\tau_1,\tau_2,c)$ in variables $a^{k-1}\in \R^{k-1},\ 
\tau_1\in \R,\tau_2\in \R,$ and $c \in \R$ such that
\beq \label{dist?}
R(a^{k-1},\tau_1,\tau_2,d_{1,P}(a_{\tau_1},a_{\tau_2}))=0.
\eneq
Moreover, $R(a^{k-1},\tau,\tau,0)\equiv 0$.
\elm
{\it{Proof}} \ \ Recall that $\{*(dP)_j(x)\}_{j=1}^k$
form an orthogonal basis in the orthogonal complement to the 
tangent plane to the level set $P^{-1}(P(x))$ at the point $x$
(see (\ref{schm})). Let us make several remarks about the rational 
function $R_{P}(x,y)$ defined by (\ref{scalar}).

1. If $a_\tau=(a^{k-1},\tau) \in K_{\m,\dt}$ is a regular
value for the map $(P_1,\dots, P_{k}):\R^N \to \R^k$ for some
$\tau$, then $a^{k-1}\in \R^{k-1}$ is a regular value for the map
$(P_1,\dots, P_{k-1}):\R^N \to \R^{k-1}$;

2. If $a^{k-1} \in \R^{k-1}$ is a regular value for the map 
$(P_1,\dots, P_{k-1}):B^N \to \R^{k-1}$, then there is
a positive $\eps(a^{k-1})>0$ such that for each point
$x \in (P_1,\dots, P_{k-1})^{-1}(a^{k-1})$ and each 
$1\leq j\leq k-1$
\beq
(*(dP)_j(x),*(dP)_j(x))>\eps(a^{k-1})>0.
\eneq
This follows from compactness of $B^N$ and regularity of 
the value $a^{k-1}$;

3. Since we consider only those $\tau$ that $a_\tau$ belongs to
the $({\m,\dt})$-cone $K_{\m,\dt}$ of regular values of $P$ 
there exists a positive constant $\eps(a^{k-1},\tau)>0$ such that 
for each  $x \in P^{-1}(a_\tau)$  
\beq
(*(dP)_k(x),*(dP)_k(x))>\eps(a^{k-1},\tau).
\eneq
This shows that $Q(x,y)$ restricted to 
$P^{-1}(K_{\m,\dt})\times P^{-1}(K_{\m,\dt})$ is a smooth function
of $x$ and $y$.

Consider irreducible representation of the rational function 
$Q(x,y)$ as a ration of two polynomials 
$Q(x,y)=\frac{T(x,y)}{S(x,y)}$. Because of remarks 2 and 3
$S(x,y) \neq 0$ for each pair 
$(x,y)\in P^{-1}(K_{\m,\dt})\times P^{-1}(K_{\m,\dt})$.
 
Now notice that we deal with smooth objects: 
smooth level sets $P^{-1}(a_\tau)$ and the smooth function
$Q_P(x,y)$.  Notice that 
$d^0_{1,P}(x,a_{\tau_2})$ is an extremal value of
the function $Q_P(x,y)$ provided that $P(y)=a_{\tau_2}$.
Similarly,  $d_{1,P}(a_{\tau_1},a_{\tau_2})$ is an extremal value 
of the function $d^0_{1,P}(x,a_{\tau_2})$ provided that 
$P(y)=a_{\tau_1}$. To find all extremal values of a smooth
function on a smooth manifold one can use  the Lagrange 
multipliers method. We prove that functions $d_{1,P}$ and $d^0_{1,P}$ 
are algebraic functions. 

The key point of the Lagrange multipliers method is that at an
extremal point of $Q_P(x,y)$ under the condition $P(y)=a_{\tau_2}$
the gradient $\nabla_yQ_P(x,y)$ can be expressed as a linear 
combination of gradients $\nabla P_1(y),\dots ,\nabla P_k(y)$,
and $\nabla r(y).$ The gradient of $Q_P(x,y)$ has the form
$$
\nabla_yQ_P(x,y)=\nabla_y\left(\frac {T_P(x,y)}{S_P(x,y)} \right)=
\left( S_P(x,y) \nabla_yT_P(x,y)- T_P(x,y) \nabla_yS_P(x,y)\right) 
S_P^{-2}(x,y).
$$

Since $S|_{P^{-1}(K_{\m,\dt})\times P^{-1}(K_{\m,\dt})}
\neq 0$ we can rewrite the Lagrange system in the following 
form
\begin{equation}\label{lagr}
\begin{cases} 
S_P(x,y) \nabla_yT_P(x,y)- T_P(x,y) \nabla_yS_P(x,y)+\\
+S_P^2(x)\left[\sum_{j=1}^k \lambda_j\nabla P_{j}(x)-
\lb_{k+1} \nabla r(y) \right]=0,\\
P(y)-a_{\tau_2} \\
T_P(x,y)-cS_P(x,y)=0, \\
\lb_{k+1} r(y)=0. 
\end{cases}
\end{equation}
{\it{Important that all equations are polynomial and we can apply
elimination theory!}}
Notice that the last equation is responsible for an extremal
point $y$ which might belong to the boundary $\partial B^N$.
If a critical value belongs to the boundary, then $r(y)=0$ and 
$\lb_{k+1} \nabla r(y)$ is not zero and 
the gradient $\nabla_y Q_P(x,y)$ should be expressed as a linear 
combination of $k+1$ vectors $\nabla P_1(y),\dots ,\nabla P_k(y)$,
and $\nabla r(y).$ If a critical value does not belong to the boundary, 
i.e. $r(y)\neq 0$, then $\lb_{k+1}=0$ and $\lb_{k+1} \nabla r(y)=0$. 

Complexify the system (\ref{lagr}), i.e. consider the system 
(\ref{lagr}) for 
$$
(x,y,\lb,a^{k-1},\tau_2,c)\in \C^N\times \C^N \times \C^{k+1} 
\times \C^{k-1} \times \C \times \C.
$$
It defines a constructible set, denoted by $V$,
in $\C^{2N+2k+2}$. Let us eliminate variables 
$\{\lb_j\}_{j=1}^{k+1},\ \{y_i\}_{i=1}^N$ by projecting 
$V$ along the corresponding $(N+k+1)$-dimensional $(\lb,y)$-plane.
The result of projection is a constructible set $W$ in
the space $(x,a^{k-1},\tau_2,c)\in \C^{N+k+1}$.
By the construction a point $(x,a^{k-1},\tau_2,c)$
belongs to $W$ if some value of $y$ the following conditions 
hold:\  
$P(y)=a_{\tau_2}$, $Q_P(x,y)=c$, and $y$ is the critical point 
of $Q_P$ restricted to $P^{-1}(a_{\tau_2})$.

The constructible set $W$ has dimension $N+k$. Indeed,
consider a polynomial function $\rho_x:P^{-1}(a_{\tau_2})\to \R$
defined  by $\rho_x(y)=Q_P(x,y)$. By Sard's lemma
for algebraic sets {\cite{Mu}} critical values of $\rho_x$ 
form an algebraic set of positive codimension in $\R$.
Therefore, the set of critical values consists of a finite number
of points and a finite number of possible $c$ so that
 $(x,a^{k-1},\tau_2,c)\in W$. Thus, $\dim W$ equals dimension
of $(x,a^{k-1},\tau_2,)$-plane, i.e. $\dim W=N+k$.

Since $W$ is constructible and has codimension $1$ there is
a non zero polynomial $\~ R(x,a_{\tau_2},c)$
such that $W \subseteq \~ R^{-1}(0)$. By the 
definition (\ref{dist2}) of $d^0_{1,P}(x,a_{\tau_2})$
and by the construction 
\beq\label{implicat}
\~ R(x,a_{\tau_1},a_{\tau_2},d^0_{1,P}(x,a_{\tau_2}))\equiv 0.
\eneq

In order to prove that the function $d_{1,P}(a_{\tau_1},a_{\tau_2})$
defined by (\ref{dist2}) is also algebraic, calculate 
critical values of  $d^0_{1,P}(x,a_{\tau_2})$, provided $P(x)=a_{\tau_1}$.
By the implicit function theorem 
the gradient $\nabla_x d^0_{1,P}(x,a_{\tau_2})$ can be 
expressed in terms of partial derivatives of  $\~R(x,a_{\tau_2},c)$
by the following way
\beq
\partial_{x_j} d^0_{1,P}(x,a)=
\partial_c \~R(x,a,c)
\left(\partial_{x_j} \~R(x,a,c)\right)^{-1},
\eneq\ 
for $(a,c)=(a_{\tau_2},d^0_{1,P}(x,a_{\tau_2}))$ and 
provided that 
$\partial_{x_j}\~R(x,a_{\tau_2},d^0_{1,P}(x,a_{\tau_2}))\neq 0$ for 
all $1\leq j \leq N$.
Fix $a=a_{\tau_2}$ and consider $x$ outside of
the union of algebraic sets 
$\Cal B=\cup_{j=1}^N\{x: \partial_{x_j}\~R(x,a,c)=0\}$.
Then $d^0_{1,P}(x,a)$ is a smooth function in $x$.
Application of the Lagrange multipliers method 
shows that at an extremal point of the function $d^0_{1,P}(x, a)$,
provided $P(x)=a_{\tau_1}$, the gradient 
$\nabla_x d^0_{1,P}(x,a)$ can  be represented as a linear combination
$\nabla_x d^0_{1,P}(x,a)=\sum_{j=1}^k \lb_j \nabla P_j(x)-
\lb_{k+1} \nabla r(x)$. Plugging in the expression 
for $\nabla_x d^0_{1,P}(x,a)$ in terms of 
$\partial_c \~ R(x,a_{\tau_2},c)$ and 
$\partial_{x_j} \~ R(x,a_{\tau_2},c)$ for $j=1,\dots ,N$ 
we can present a Lagrange multiplier system in the 
following form
\beq \label{lagr1}
\begin{cases} 
\partial_{x_j} \~ R(x,a_{\tau_2},c) 
\partial_c \~ R(x,a_{\tau_2},c)=\\
= \left[\sum_{j=1}^k \lambda_j\nabla P_{j}(x)-
\lb_{k+1} \nabla r_1(x) \right] 
\left(\partial_{x_j} \~ R(x,a_{\tau_2},c)\right)^2\\ 
\~ R(x,a_{\tau_2},c)=0,\\
P(x)=a_{\tau_1},\\
\lb_{k+1} r(x)=0.
\end{cases}
\eneq
Again the system (\ref{lagr1}) consists of only polynomial equations and
we can apply elimination theory.
Consider this system  for 
$$
(x,\lb,a^{k-1},\tau_1,\tau_2,c)\in \C^N\times  \C^{k+1} 
\times \C^{k-1} \times \C \times \C \times \C.
$$ 
It defines a constructible set, denoted by $V_1$,
in $\C^{N+2k+3}$. Let us eliminate variables 
$\{\lb_j\}_{j=1}^{k+1},\ \{x_i\}_{i=1}^N$ by projecting 
$V$ along the corresponding $(N+k+1)$-dimensional $(\lb,x)$-plane.
The result of projection is a constructible set $W_1$ in
the space $(a^{k-1},\tau_1,\tau_2,c)\in \C^{k+2}$.

Similarly to the arguments for the constructible set $W$
one can show that $W_1$ has dimension $k$. 
Since $W_1$ is constructible and has codimension $1$ there is
a nonzero polynomial $R(a^{k-1},\tau_1,\tau_2,c)$
such that $W_1 \subseteq R^{-1}(0)$. By the 
definition (\ref{dist2}) of $d_{1,P}(a_{\tau_1},a_{\tau_2})$
and by the construction 
\beq\label{implic}
R(a_{\tau_1},a_{\tau_2},d_{1,P}(a_{\tau_1},a_{\tau_2}))\equiv 0.
\eneq

By the construction if $a_{\tau_1}=a_{\tau_2}$,
then $R(a_{\tau_1},a_{\tau_1},0)\equiv 0$,
because in this case both level sets are the same and
$C^1$-distance between them must equal zero. This completes 
the proof of Lemma \ref{polynom}.

\blm With the notations above there exists a refinement 
$(\m',\dt')$-cone $K_{\m',\dt'}\subset K_{\m,\dt}$ such that 
for any pair of points $a_{\tau_1}=(a^{k-1},\tau_1)$ and 
$a_{\tau_2}=(a^{k-1},\tau_2)$ from $K_{\m',\dt'}$ 
\beq
D^1_P(a_{\tau_1},a_{\tau_2})\leq (a_1 \dots a_1{k-1})^2,
\eneq
where $a^{k-1}=(a_1, \dots ,a_{k-1}) \in \R^{k-1}$.
\elm

{\it{Proof}}\ \ It follows from Lemma \ref{polynom} that 
there is a polynomial $R(a_{\tau_1},a_{\tau_2},c)$ such that
$R(a_{\tau_1},a_{\tau_2},d_{1,P}(a_{\tau_1},a_{\tau_2}))\equiv 0$
and $R(a,\tau,\tau,0)\equiv 0$.

For $a_{\tau_1},a_{\tau_2}$ belonging to the $(\m,\dt)$-cone
$K_{\m,\dt}$ of regular values of $P$ 
the function $d_{1,P}(a_{\tau_1},a_{\tau_2})$ depends continuously
on $\tau_1$ and $\tau_2$. Let us rewrite $R(a_{\tau_1},a_{\tau_2},c)$
in the form $R(a^{k-1},\tau_1,\tau_2,c)$.
Recall that in our notations $a_{\tau}=(a^{k-1},\tau)$.

Suppose for determiness that $\tau_1>\tau_2$. Notice that
each sufficiently small positive root $c_j(a^{k-1},\tau_1,\tau_2)$
is increasing in $\tau_1$ and decreasing in $\tau_2$ in a 
neighborhood of $(a^{k-1},\tau_1,\tau_2)=0$. 
Therefore, $c_j(a^{k-1},\tau_1,0)>c_j(a^{k-1},\tau_1,\tau_2)$

Denote $R(a^{k-1},\tau,c)=R(a^{k-1},\tau,0,c)$.
Let us show that for some sufficiently large positive 
integers $m'$ and $m''$ 
if $0<\tau<(a^{k-1}_1 \dots a^{k-1}_{k-1} c)^{m'}$
and $0<c<(a^{k-1}_1 \dots a^{k-1}_{k-1})^{m''}$ the following
decomposition holds:
Put $c=\rho (a^{k-1}_1 \dots a^{k-1}_{k-1})^{m''}$ and
$\tau=\lb(a^{k-1}_1 \dots a^{k-1}_{k-1} c)^{m'}$ with
$\rho,\lb \in (0,1)$. Then
\beq 
R(a^{k-1},\tau,c)=r_0(a^{k-1})(1+q_1(a^{k-1},\rho))(1+q_2(a^{k-1},\rho,\lb)),
\eneq
where $r_0,\ q_1,\ q_2$ are polynomials in their variables.
Indeed, apply the same arguments as we used to prove (\ref{prod}) to
the polynomial
$$
R(a^{k-1},\tau,c)=\sum_{j=1}^{\deg R}R_j(a^{k-1},c)\tau^j+R_0(a^{k-1},c).
$$
Then apply the same arguments to 
$$
R_0(a^{k-1},c)=\sum_{j=1}^{\deg R_0}r_j(a^{k-1})c^j+r_0(a^{k-1},c).
$$
Notice that $R(a^{k-1},0,0)\equiv 0$ implies that 
$q_1(0,\rho)\equiv q_2(0,\rho,\lb)\equiv 0$.
Therefore, for a sufficiently small $\dt'>0$ and any 
$(a^{k-1},\dt')\in K_{\m,\dt'}$ polynomials 
$q_1(a^{k-1},\rho)$ and  $q_2(a^{k-1},\rho,\lb)$ are sufficiently small
and $R(a^{k-1},\tau,c)$ equals $0$ if  and only if
$r_0(a^{k-1})$ equals $0$.

By Lemma {cone} there is a refinement $(\m',\dt')$-cone 
$K_{\m',\dt'} \subset K_{\m,\dt'}$ such that $r_0(a^{k-1})$
does not vanish on $K_{\m',\dt'}$.

Now put $m_k'=m' m''$. As we have just shown all sufficiently
small positive roots $c_j(a^{k-1},\tau,0)<\tau^{1/m'}
(a^{k-1}_1 \dots a^{k-1}_{k-1})$ provided that 
$$
\tau^{1/m'}(a^{k-1}_1 \dots a^{k-1}_{k-1})<
(a^{k-1}_1 \dots a^{k-1}_{k-1})^{m''}.
$$
This condition is satisfied for any 
$0<\tau<(a^{k-1}_1 \dots a^{k-1}_{k-1})^{m' m''}$.
This shows that all sufficiently small positive roots
$$
c_j(a^{k-1},\tau,0)<(a^{k-1}_1 \dots a^{k-1}_{k-1})^{m''+1}<
(a^{k-1}_1 \dots a^{k-1}_{k-1})^2.
$$
This completes the proof of the Lemma.

Let us complete the proof of part c) Theorem \ref{existence1}
by the following inductive arguments.

Consider a sequence of positive integers $m_2,\dots , m_k$.
Let $\m=(1,m_2,\dots , m_k)$ and $\dt>0$.
Define a sequence of polynomials associated to this sequence, 
defined by their coordinate functions:
\beq 
\begin{aligned}
P_j^0=P_j-(P_1\dots P_{j-1})^{m_j},\ j=2,\dots ,k\\
P^s=(P_1^0,\dots,P_s^0,P_{s+1},\dots ,P_{k-s}),\ s=2,\dots ,k
\end{aligned}
\eneq
Define the restriction of the $(\m,\dt)$-cone $K_{\m,\dt}$ to the 
$s$-dimensional plane, denoted by $K^s_{\m,\dt}$, generated by the 
first $k$-coordinates by the following way:
\beq 
\begin{aligned}
K^s_{\m,\dt}=\{a^s=(a_1,\dots, a_s)\in \R^s: \ 
0<a_1<\dt, \\
0 < |a_{j+1}| < |a_1 \dots a_j|^{m_{j+1}}\ 
\text{for}\ j=1, \dots, s-1\}.
\end{aligned}
\eneq 

It is shown above that there is an $(\m,\dt)$-cone $K_{\m,\dt}$ 
such that any point $(0,a^{k-1})\in \R \times K^{k-1}_{\m,\dt}$ 
is a regular point for the polynomial $P^{k-1}$.
Therefore, one can apply Lemmas \ref{cone}, \ref{polynom},
and \ref{dist?} and show that there is refinement of
$K^{k-1}_{\m,\dt}$, denoted the same, such that for any two 
points $a^{k-1}_{\tau_1}=(0,a^{k-2},\tau_1)$ and 
$a^{k-1}_{\tau_1}(0,a_{k-2},\tau_2)$ from $\R \times K^{k-1}_{\m,\dt}$ 
\beq
D^1_{P^1}(a^{k-1}_{\tau_1},a^{k-1}_{\tau_2})\leq (a_1 \dots a_{k-2})^2.
\eneq
By induction one can show that there is a refinement
an $(\m,\dt)$-cone $K_{\m,\dt}$ such that for any  two 
points $a^{k-s}_{\tau_1}=(0,a^{k-s-1},\tau_1)$ and 
$a^{k-s}_{\tau_1}(0,a_{k-s-1},\tau_2)$ from the restriction cone
$K^{k-s}_{\m,\dt}$
such that 
\beq
D^1_{P^s}(a^{k-s}_{\tau_1},a^{k-s}_{\tau_2})\leq (a_1 \dots a_{k-s-1})^2.
\eneq
Notice that for any $1\leq s\leq k$ level sets of the  polynomial 
$P^s$ correspond to level sets of the initial polynomial $P$.
Combining this with all estimates for 
$D^1_{P^s}(a^{k-s}_{\tau_1},a^{k-s}_{\tau_2})$
and the triangle-like inequality (\ref{pstriang}) 
one can show that part c) of Theorem \ref{existence1} holds true.
Part a) of Theorem \ref{existence1} follows from part c)
because $Q_P(x,y)\geq \|x-y\|^2$, which implies that
$D^1_P(a,b)\geq D^0_P(a,b)$ for any pair $a,b \in K_{\m,\dt}$. 

This completes the proof of Theorem \ref{existence1}.

\end{document}